\def\DATE{\relax}
\magnification=1100
\baselineskip=12.72pt
\voffset=.75in
\hoffset=.9in
\hsize=4.1in
\newdimen\hsizeGlobal
\hsizeGlobal=4.1in%
\vsize=7.05in
\parindent=.166666in
\pretolerance=500 \tolerance=1000 \brokenpenalty=5000

\footline={\vbox{\hsize=\hsizeGlobal\hfill{\rm\the\pageno}\hfill\llap{\sevenrm\DATE}}\hss}

\long\def\note#1{%
  \hfuzz=50pt%
  \vadjust{%
    \setbox1=\vtop{%
      \hsize 3cm\parindent=0pt\eightpoints\baselineskip=9pt%
      \rightskip=4mm plus 4mm\raggedright#1\hss%
      }%
    \hbox{\kern-4cm\smash{\box1}\hss\par}%
    }%
  \hfuzz=0pt
  }
\def\note#1{\relax}

\def\anote#1#2#3{\smash{\kern#1in{\raise#2in\hbox{#3}}}%
  \nointerlineskip}     

\newcount\equanumber
\equanumber=0
\newcount\sectionnumber
\sectionnumber=0
\newcount\subsectionnumber
\subsectionnumber=0
\newcount\snumber  
\snumber=0

\def\section#1{%
  \subsectionnumber=0%
  \snumber=0%
  \equanumber=0%
  \advance\sectionnumber by 1%
  \noindent{\bf \the\sectionnumber .~#1{}.~}%
}%
\def\subsection#1{%
  \advance\subsectionnumber by 1%
  \snumber=0%
  \equanumber=0%
  \noindent{\bf \the\sectionnumber .\the\subsectionnumber .~#1.~}%
}%
\def\prevs{\the\sectionnumber .\the\subsectionnumber .\the\snumber }
\long\def\Definition#1{%
  \global\advance\snumber by 1%
  \bigskip
  \noindent{\bf Definition~\prevs .}%
  \quad{\it#1}%
}

\long\def\Corollary#1{%
  \global\advance\snumber by 1%
  \bigskip
  \noindent{\bf Corollary~\prevs .}%
  \quad{\it#1}%
}%
\long\def\Lemma#1{%
  \global\advance\snumber by 1%
  \bigskip
  \noindent{\bf Lemma~\prevs .}%
  \quad{\it#1}%
}%
\def\Proof{\noindent{\bf Proof.~}}
\long\def\Proposition#1{%
  \advance\snumber by 1%
  \bigskip
  \noindent{\bf Proposition~\prevs .}%
  \quad{\it#1}%
}%
\long\def\Remark#1{%
  \bigskip
  \noindent{\bf Remark.~}#1%
}%
\long\def\remark#1{%
  \advance\snumber by1%
  \bigskip
  \noindent{\bf Remark~\prevs .}\quad#1%
}%
\long\def\Theorem#1{%
  \advance\snumber by 1%
  \bigskip
  \noindent{\bf Theorem~\prevs .}%
  \quad{\it#1}%
}%
\long\def\Statement#1{%
  \advance\snumber by 1%
  \bigskip
  \noindent{\bf Statement~\prevs .}%
  \quad{\it#1}%
}%
\def\ifundefined#1{\expandafter\ifx\csname#1\endcsname\relax}
\def\labeldef#1{\global\expandafter\edef\csname#1\endcsname{\prevs}}
\def\labelref#1{\expandafter\csname#1\endcsname}
\def\label#1{\ifundefined{#1}\labeldef{#1}\note{$<$#1$>$}\else\labelref{#1}\fi}

\def\preveq{(\the\sectionnumber .\the\subsectionnumber .\the\equanumber)}
\def\neq{\global\advance\equanumber by 1\eqno{\preveq}}

\def\ifundefined#1{\expandafter\ifx\csname#1\endcsname\relax}

\def\equadef#1{\global\advance\equanumber by 1%
  \global\expandafter\edef\csname#1\endcsname{\preveq}%
  \setbox1=\hbox{\rm\hskip .1in[#1]}\dp1=0pt\ht1=0pt\wd1=0pt%
  \preveq\box1}
\def\equadef#1{\global\advance\equanumber by 1%
  \global\expandafter\edef\csname#1\endcsname{\preveq}%
  \preveq}

\def\equaref#1{\expandafter\csname#1\endcsname}

\def\equa#1{%
  \ifundefined{#1}%
    \equadef{#1}%
  \else\equaref{#1}\fi}

\font\eightrm=cmr8%
\font\sixrm=cmr6%

\font\eightsl=cmsl8%

\font\eightbf=cmb8%

\font\eighti=cmmi8%
\font\sixi=cmmi6%

\font\eightsy=cmsy8%
\font\sixsy=cmsy6%

\font\eightex=cmex8%
\font\sixex=cmex6%
\font\fiveex=cmex5%

\font\eightit=cmti8%

\font\eighttt=cmtt8%

\font\tenbb=msbm10%
\font\eightbb=msbm8%
\font\sevenbb=msbm7%
\font\sixbb=msbm6%
\font\fivebb=msbm5%
\newfam\bbfam  \textfont\bbfam=\tenbb  \scriptfont\bbfam=\sevenbb  \scriptscriptfont\bbfam=\fivebb%

\font\tenbbm=bbm10

\font\tencmssi=cmssi10%
\font\sevencmssi=cmssi7%
\font\fivecmssi=cmssi5%
\newfam\ssfam  \textfont\ssfam=\tencmssi  \scriptfont\ssfam=\sevencmssi  \scriptscriptfont\ssfam=\fivecmssi%

\font\tenfrak=cmfrak10%
\font\eightfrak=cmfrak8%
\font\sevenfrak=cmfrak7%
\font\sixfrak=cmfrak6%
\font\fivefrak=cmfrak5%
\newfam\frakfam  \textfont\frakfam=\tenfrak  \scriptfont\frakfam=\sevenfrak  \scriptscriptfont\frakfam=\fivefrak%
\def\frak{\fam\frakfam\tenfrak}%

\font\tenmsam=msam10%
\font\eightmsam=msam8%
\font\sevenmsam=msam7%
\font\sixmsam=msam6%
\font\fivemsam=msam5%

\def\bb{\fam\bbfam\tenbb}%

\def\hexdigit#1{\ifnum#1<10 \number#1\else%
  \ifnum#1=10 A\else\ifnum#1=11 B\else\ifnum#1=12 C\else%
  \ifnum#1=13 D\else\ifnum#1=14 E\else\ifnum#1=15 F\fi%
  \fi\fi\fi\fi\fi\fi}
\newfam\msamfam  \textfont\msamfam=\tenmsam  \scriptfont\msamfam=\sevenmsam  \scriptscriptfont\msamfam=\fivemsam%
\def\msam{\msamfam\tenmsam}%
\mathchardef\leq"3\hexdigit\msamfam 36%
\mathchardef\geq"3\hexdigit\msamfam 3E%

\font\tentt=cmtt11%
\font\seventt=cmtt9%
\textfont\ttfam=\tentt
\scriptfont7=\seventt%
\def\tt{\fam\ttfam\tentt}%

\def\eightpoints{%
\def\rm{\fam0\eightrm}%
\textfont0=\eightrm   \scriptfont0=\sixrm   \scriptscriptfont0=\fiverm%
\textfont1=\eighti    \scriptfont1=\sixi    \scriptscriptfont1=\fivei%
\textfont2=\eightsy   \scriptfont2=\sixsy   \scriptscriptfont2=\fivesy%
\textfont3=\eightex   \scriptfont3=\sixex   \scriptscriptfont3=\fiveex%
\textfont\itfam=\eightit  \def\it{\fam\itfam\eightit}%
\textfont\slfam=\eightsl  \def\sl{\fam\slfam\eightsl}%
\textfont\ttfam=\eighttt  \def\tt{\fam\ttfam\eighttt}%
\textfont\bffam=\eightbf  \def\bf{\fam\bffam\eightbf}%

\textfont\frakfam=\eightfrak  \scriptfont\frakfam=\sixfrak \scriptscriptfont\frakfam=\fivefrak  \def\frak{\fam\frakfam\eightfrak}%
\textfont\bbfam=\eightbb      \scriptfont\bbfam=\sixbb     \scriptscriptfont\bbfam=\fivebb      \def\bb{\fam\bbfam\eightbb}%
\textfont\msamfam=\eightmsam  \scriptfont\msamfam=\sixmsam \scriptscriptfont\msamfam=\fivemsam  \def\msam{\msamfam\eightmsam}

\rm%
}

\newdimen\poorBoldDim
\def\poorBold#1{\setbox1=\hbox{#1}\poorBoldDim=\wd1\wd1=0pt\copy1\kern0.25pt\copy1\kern0.25pt\box1\kern\poorBoldDim}%

\mathchardef\lsim"3\hexdigit\msamfam 2E%
\mathchardef\gsim"3\hexdigit\msamfam 26%

\def\d{\,{\rm d}}

\def\ds{\displaystyle}
\long\def\DoNotPrint#1{\relax}

\def\finetune#1{#1}
\def\fixedref#1{#1\note{fixedref$\{$#1$\}$}}

\def\limn{\lim_{n\to\infty}}
\def\limN{\lim_{N\to\infty}}
\def\liminfn{\liminf_{n\to\infty}}
\def\limsupn{\limsup_{n\to\infty}}

\def\qed{{\vrule height .9ex width .8ex depth -.1ex}}
\def\sign{\mathop{\rm sign}}
\def\ss{\scriptstyle}
\def\step#1{\medskip\noindent{\bf #1.}\ }

\def\boc{\note{{\bf BoC}\hskip-11pt\setbox1=\hbox{$\Bigg\downarrow$}%
         \dp1=0pt\ht1=0pt\ht1=0pt\leavevmode\raise -20pt\box1}}
\def\eoc{\note{{\bf EoC}\hskip-11pt\setbox1=\hbox{$\Bigg\uparrow$}%
         \dp1=0pt\ht1=0pt\ht1=0pt\leavevmode\raise 20pt\box1}}

\def\One{\hbox{\tenbbm 1}}

\def\calB{{\cal B}}
\def\calC{{\cal C}}

\def\MM{{\bb M}\kern .4pt}
\def\NN{{\bb N}\kern .5pt}

\pageno=1

\centerline{\bf SOME TAUBERIAN THEORY}
\centerline{\bf FOR THE \poorBold{$q$}-LAGRANGE INVERSION}

\bigskip
 
\centerline{Ph.\ Barbe$^{(1)}$ and W.P.\ McCormick$^{(2)}$}
\centerline{${}^{(1)}$CNRS {\sevenrm(UMR {\eightrm 8088})}, ${}^{(2)}$University
 of Georgia}

 
{\narrower
\baselineskip=9pt\parindent=0pt\eightpoints

\bigskip

{\bf Abstract.} We consider formal power series defined through the
functional $q$-equation of the $q$-Lagrange inversion. Under some assumptions,
we obtain the asymptotic behavior of the coefficients of these power 
series. As a by-product, we show that, via the $1/q$-Borel transform, 
the $q$-Lagrange inversion formula provides an interpolation between the usual
Lagrange inversion, $q=1$, and the probabilistic theory of renewal
sequences, $q\to 0$. We also discuss some new solutions of the $q$-Lagrange
inversion equation which do not vanish at $0$.

\bigskip


\bigskip
 

}

\bigskip

\def\prevs{\the\sectionnumber.\the\snumber }
\def\preveq{(\the\sectionnumber.\the\equanumber)}

\section{Introduction}
In their $q$-extension of the Lagrange inversion formula, Andrews (1975),
Gessel (1980),
Garsia (1981) and Krattenthaler (1988) among others consider the following
problem: let $f(z)=\sum_{n\geq 1} f_n z^n$ be a formal power series with
$f_1\not=0$ and define 
the linear operator $U_{f,q}$ on formal power series by extending linearly
$U_{f,q}z^k=f(z)f(qz)\cdots f(q^{k-1}z)$, $k\in\NN$, with $U_{f,q}z^0=1$; what is
the inverse of $U_{f,q}$? They show that if one
defines the formal power series $g(z)$ by
$$
  \sum_{n\geq 1} f_n g(z)g(z/q)\cdots g(z/q^{n-1})=z \, ,
  \eqno{\equa{rightInverseDef}}
$$
and the requirement that the constant term of $g$ vanishes,
then the inverse of $U_{f,q}$ is $U_{g,1/q}$. While this result is of 
considerable interest in combinatorics and in the theory of $q$-series,
there is no simple way to determine the function $g$ besides calculating
its coefficient recursively from \rightInverseDef\ or using Garsia and
Haiman's (1996) formula. In any cases, the power series $g$ remains 
rather mysterious and very little is known about it.

When $q$ is $1$, \rightInverseDef\ is $f\bigl(g(z)\bigr)=z$, and, even if
one may not be able to calculate $g$ explicitly, one can often relate the
asymptotic behavior of its coefficients to some accessible features of $f$. 
This viewpoint is very useful in combinatorics, as stressed for instance by 
Flajolet and Sedgewick (2009; \S VI.7, \S VII.3).

The purpose of this paper is to determine the asymptotic behavior of
the coefficients of $g$ involved in \rightInverseDef\ in terms of the
function $f$, under suitable conditions on $f$. These conditions may
force $g$ to be a divergent series, which precludes the use of Cauchy
integral formula and traditional complex analytic methods to derive
the asymptotic behavior of the coefficients of $g$. Thus, our paper
may also be viewed as a contribution to the asymptotic methods for
divergent series.

Some information on the asymptotic behavior of the coefficients of $g$ can 
be deduced from more general results by
Zhang (1998) and Cano and Fortuny Ayuso (2012) on $q$-analytic equations; 
but these results are in 
terms of $q$-Gevrey orders and, while obtained under rather general 
assumptions, they deliver only rather crude bounds for our specific problem.
More refined results on $q$-algebraic equation were obtained by Barbe and
McCormick (2013), which can be used when only a finite number of
coefficients $f_n$ do not vanish. On the contrary, in the current paper
we are interested in situations where the $f_n$ do not vanish ultimately, and 
the structure of the $q$-Lagrange inversion will allow us to derive, in this
context, more precise results than those in Barbe and McCormick (2013).
For polynomials arising from counting certain lattice paths by area,
Drake (2009) obtained asymptotic results on coefficients of some power 
series satisfying a functional relation rather similar though not equivalent
to \rightInverseDef. In contrast, our approach is mostly analytical.

We will also see that a change of function allows one to 
linearize \rightInverseDef, and relate its study to that of linear 
$q$-difference equation.

Following Garsia (1981), if $f$ and $g$ are related through \rightInverseDef, 
we call $g$ the right inverse of $f$.

The paper is organized as follows. In section 2, we state our assumptions and
main results when $0<q<1$. Section 3 discusses the $q$-Catalan numbers, which provide a toy
example for the arguments used in the proof of our main result. Section 4 is
devoted to some examples which illustrate the diversity of asymptotic
behavior that equation \rightInverseDef\ yields, when $0<q<1$. In section 5, 
we provide some results when $q>1$. In section 6, we discuss the
connection between \rightInverseDef\ and the theory of linear $q$-difference
equations and provide explicit formulas for solutions of \rightInverseDef\
which do not vanish at $0$. Section 7 is devoted to the proof of the main
result of this paper, and the remaining sections contain proofs of further
results stated in sections 2 and 5.

\bigskip
\finetune{\vfill\eject}
\section{Results for \poorBold{$0<q<1$}}
For any real number $c$, consider the composition operator $M_cg(z)=g(cz)$.
Since $U_{cf,q}z^k=U_{f,q}M_c z^k$, we see that if $U_{f,q}^{-1}=U_{g,1/q}$, then
for any $c$ nonnegative, 
$$
  U_{cf,q}^{-1}
  =M_{1/c}U_{g,1/q}
  =U_{M_{1/c}g,1/q} \, .
$$
Therefore, up to considering $f/f_1$, we can
assume without any loss of generality that $f_1=1$. When $f$ is
linear, this means $f(z)=z$, and, in this case, $U_{f,q}z^k=q^{k\choose 2} z^k$
has inverse $U_{f,q}^{-1}z^k=q^{-{k\choose 2}} z^k$, so that $g(z)=z$ as well.
Once the linear case is ruled out, perhaps the next simplest one is the 
quadratic
one, where $f(z)=z(1-z)$. The importance of this case stems from its bearing
to the theory of $q$-series, for which, with the traditional notation
$(z;q)_n=\prod_{0\leq j<n} (1-q^jz)$, we 
have $U_{f,q}z^k=q^{k\choose 2} z^k(z;q)_k$. In this case, neither the 
function $g$ nor its coefficients $(g_n)$ are known in a closed form, 
illustrating the complex relationship between $f$ and $g$. 

Whenever we have a power series $f(z)$, we write $[z^i]f$ for its $i$-th
coefficient, that is for $f_i$. Recall that a function $h$ is regularly
varying at infinity of index $\rho$ if for any $\lambda$ positive,
$$
  \lim_{t\to\infty} h(\lambda t)/h(t)=\lambda^\rho \, .
$$
We refer to Bingham, Goldie and Teugels (1989) for the theory of regularly 
varying functions. The notion of regularly varying function is extended
to that of regularly varying sequence by considering that a sequence is a 
function constant on intervals $[n,n+1)$, $n\in\NN$.

To state our result, recall that we assume without loss of generality that
$f_1=1$. We then define the power series $\phi$ by
$$
  f(z)=z\bigl(1-\phi(z)\bigr) \, .
  \eqno{\equa{phiDef}}
$$
In particular, $\phi(0)=0$. Throughout this paper, we set $\phi_n=[z^n]\phi$.

Our first result sets the scene for the main contribution of this paper.
In order to keep the exposition flowing we defer its proof to the last section
of this paper.

\Proposition{\label{gDivergent}
  Assume that $\phi\not=0$, and that either $(\phi_n)$ 
  or $\bigl( (-1)^{n+1}\phi_n\bigr)$ is a nonnegative sequence. Then, any
  formal solution of \rightInverseDef\ is a divergent series.
}

\bigskip

If \rightInverseDef\ holds, the coefficients $(g_n)$ of $g$ are power series
in $1/q$.  Our next result show that they are in fact polynomials in $1/q$
and that their leading term is related to the $n$-th coefficient of $z/f(z)$.

Whenever we have a power series $h(q)$, we write $[q^j]h$ for the 
coefficient of $q^j$ in $h$.

\Proposition{\label{gnPolynomial}
  Let $g$ be a solution of \rightInverseDef\ (with $g_0=0$). For any $n$ the 
  coefficient $g_n=[z^n]g$ is a polynomial in $1/q$, of degree
  at most ${n\choose 2}$. Furthermore, 
  $$
    [1/q^{n\choose 2}]g_n = [z^{n-1}]{1\over 1-\phi} \, ,
  $$
  and, if $(\phi_n)$ is a nonnegative sequence, the coefficients
  of the polynomial $g_n$ (in $1/q$) are nonnegative. 
}

\bigskip

It follows from Proposition \gnPolynomial\ that $q^{n\choose 2} g_n$ is a 
polynomial in $q$. If the $\phi_n$ are positive, the leading term of this 
polynomial is $q^{n\choose 2} [z^n]1/\bigl(1-\phi(z)\bigr)$. Whenever 
the equation $\phi(z)=1$ has a positive solution $\zeta$, the power 
series $1/\bigl(1-\phi(z)\bigr)$ has radius of convergence $1/\zeta$, 
and therefore $g_n$ is expected to grow like $q^{-{n\choose 2}}\zeta^{n(1+o(1))}$ 
as $n$
tends to infinity. To describe more precisely the asymptotic behavior
of the polynomial $q^{n\choose 2}g_n$ as $n$ tends
to infinity, we need the following definition. With a slight abuse of 
language, it recalls 
three modes of convergence which will be useful to state our result.

\Definition{\label{convergenceMeaning}
  Let $\bigl(a_n(q)\bigr)_{n\geq 0}$ be a sequence of polynomials in $q$. We say
  that this sequence converges to a power series $L(q)$ for any $q$ in a 
  set $A$

  \medskip
  \noindent (i) in limit infimum if for any $q$ in $A$,
  $$
    \liminf_{n\to\infty} a_n(q)=L(q) \, ;
  $$

  \noindent (ii) except on a set of integers of density $0$ if for any positive
  $\delta$ and any $q$ in $A$,
  $$
    \limn {1\over n} \sum_{1\leq i\leq n} \One\{\, |a_i(q)-L(q)|\geq\delta\,\}
    = 0 \, ;
  $$

  \noindent (iii) coefficientwise, if for any fixed nonnegative integer $j$,
  $$
    \limn [q^j]a_n=[q^j]L \, .
  $$
}

Since $\phi(0)=0$, the product
$\prod_{j\geq 0}\bigl(1-\phi(q^jz)\bigr)$ is convergent whenever $|q|<1$
and $z$ is in the interior of the disk of convergence of $\phi$.
We can now state our main result, whose assumptions imply that $\phi$ has 
a positive radius of convergence.

\Theorem{\label{mainTh}
  Assume that 

  \smallskip
  \noindent (i) $\phi_i\geq 0$ for any $i\geq 1$, 
  
  \smallskip
  \noindent (ii) $\zeta=\min\{\, x>0\,:\, \phi(x)=1\,\}$ exists;

  \smallskip
  \noindent (iii) for any nonnegative integer $j$, $\limn \phi_{j+n}/\phi_n\in
  [\,0,\infty)$.

  \smallskip
  \noindent (iv) the sequence $\Bigl([z^n]{\ds 1\over\ds 1-\phi(\zeta z)}\Bigr)_{n\geq 0}$
  is regularly varying with index $\rho-1$ in $(-1,0\,]$.

  \smallskip
  \noindent
  Then, for any $q$ in $(0,1)$, the following convergence holds in the three
  modes of convergence, limit infimum, except on a set of integers of 
  density $0$ and coefficientwise:
  $$
    \zeta^n n
    \biggl(1-\phi\Bigl(\zeta\Bigl(1-{1\over n}\Bigr)\Bigr)\biggr) \, 
    q^{n\choose 2} g_n 
    \to {\zeta\over\Gamma(\rho) \prod_{j\geq 1} \bigl( 1-\phi(\zeta q^j)\bigr)}
    \, .
  $$
}%
%

Before commenting on the meaning of the conclusion of Theorem \mainTh, we
comment on assumption (iv). Set $\tilde\phi(z)=\phi(\zeta z)$. Karamata's
theorem for power series (see Bingham, Goldie and Teugels, 1989, Corollary
1.7.3) implies that under (iv), 
$$
  [z^n]{1\over 1-\tilde\phi} 
  \sim {1\over \Gamma(\rho)n\bigl( 1-\tilde\phi(1-1/n)\bigr) }
  \eqno{\equa{mainThA}}
$$
as $n$ tends to infinity. Thus, the sequence 
\hbox{$n\bigl(1-\tilde\phi(1-1/n)\bigr)$, $n\geq 1$}, is regularly varying 
of index
$1-\rho$. Since $\tilde\phi$ is monotone under assumption (i) of Theorem 
\mainTh, this implies that the function $1-\tilde\phi(1-1/t)$ is regularly
varying of index $-\rho$ (note that Karamata's theorem, as stated in Bingham,
Goldie and Teugels, 1989, gives the stronger result that assumption (iv) of 
Theorem \mainTh\ is equivalent to \mainThA\ 
if $[z^n]\bigl(1/(1-\tilde\phi)\bigr)$, $n\geq 1$,
is a monotone sequence; however, it can be seen from the proof in Bingham,
Goldie and Teugels, 1989, that one does not need this monotonicity assumption 
in order to derive \mainThA\ from assumption (iv)).

We can now see that it is not possible to have $\rho> 1$. Indeed, if $\rho>1$,
assumption (i) of Theorem \mainTh\ implies that for any $n$ nonnegative and
any $s$ in $(0,1)$,
$$\eqalign{
  1-\tilde\phi(1-s)
  &{}=\sum_{k\geq 0} \phi_k \zeta^k \bigl(1-(1-s)^k\bigr) \cr
  &{}\geq \sum_{0\leq k\leq n} \phi_k \zeta^k \bigl( 1-(1-s)^k\bigr) \, . \cr
  }
$$
Dividing this inequality by $s$ and taking limit infimum as $s$ tends 
to $0+$, we would obtain
$$
  0\geq \sum_{0\leq k\leq n} k\phi_k \zeta^k \, .
$$
Under assumption (i), this forces $\phi_k=0$ for any $k\geq 1$. 
Since $\phi(0)=0$ under \phiDef, $\phi$
is then the constant function $0$, and assumption (ii) cannot hold. 
Consequently, we must have $\rho\leq 1$.

To summarize, assumption (iii) means practically that the 
function $1-\phi\bigl( \zeta(1-1/t)\bigr)$ is a well behaved regularly 
varying function; its index of regular variation is then some negative
number $-\rho$; this index is at least $-1$ and when the function $\phi$
is differentiable at $\zeta$, this index is $-1$.

Qualitatively, the conclusion of Theorem \mainTh\ asserts that
for such a function $\phi$,
the asymptotic behavior of $(g_n)$ is driven by a super-exponential term which
depends only on $q$, namely, $1/q^{n\choose 2}$, corrected by a term of 
algebraic decay or growth, $1/\zeta^n$, which depends on the location of the 
first zero of $f$,  and another term of algebraic 
growth, $f\bigl(\zeta(1-1/n)\bigr)$, which depends on the singularity 
of $f$ at this first positive zero. The constant term,
involves both the zero $\zeta$ and the complete function $f$ through the
product  of the nonlocal terms $1-\phi(q^j\zeta)$. 

If $f$ is 
differentiable at $\zeta$, then $\rho=1$ and 
$$
  n\biggl( 1-\phi\Bigl(\zeta\Bigl(1-{1\over n}\Bigr)\Bigr)\biggr)
  \sim \zeta \phi'(\zeta)
$$ 
as $n$ tends to infinity. This yields a simpler expression for the 
conclusion of the theorem.

\bigskip

Theorem \mainTh\ is not entirely satistifactory, for
its conclusion leaves open the possibility that the sequence involved in this
conclusion does not 
converge because of an exceptional set of integers. We will comment further on
this point at the end of this section. When the sequence $(\phi_i)$ tends to $0$
at a superexponential rate, a more
definitive conclusion can be given, which is the purpose of the following
result.

\Proposition{\label{mainThPolynomial}
  Assume that

  \medskip

  \noindent
  (i) $(\phi_n)$ is a nonnegative sequence;

  \medskip
  \noindent
  (ii) $\limsup_{n\to\infty} {\ds\log\phi_n\over\ds n\log n}<-2$;

  \medskip
  \noindent
  (iii) $\phi(z)=1-\bigl(1-(z/\zeta)\bigr)\theta(z)$ where $\theta$ is a 
  meromorphic function having no root in the closed disk of radius $\zeta$.

  \medskip
  \noindent
  Then
  $
    \ds g_n
    \sim {1\over q^{n\choose 2} \zeta^n \phi'(\zeta)(q;q)_\infty \prod_{j\geq 1}
                 \theta(\zeta q^j)}
  $
  as $n$ tends to infinity.
}

\bigskip

Proposition \mainThPolynomial\ covers the case where $\phi$ is a 
polynomial with nonnegative coefficients that sum to $1$. 
Note that \phiDef\ and the requirement that $f_1=1$ force $\phi(0)=0$, so that
assumption (iii) of Proposition \mainThPolynomial\ yields $\theta(0)=1$.
Then, all the products involved in the conclusion of Proposition 
\mainThPolynomial\ are well defined.

B\'ezivin (1992) and Ramis (1992) in their study of linear $q$-difference 
equa\-tions introduced the definition of $q$-Gevrey order with $q>1$, which 
we adapt here to the range $0<q<1$ as follows: a 
formal power series $\sum g_n z^n$ is
of $1/q$-Gevrey order $s$ if there exists a positive constant $A$ such 
that $g_n=O(A^nq^{-sn^2/2})$ as $n$ tends to infinity.
Proposition \mainThPolynomial\ implies that if its assumptions are satisfied,
then $g$ is of $1/q$-Gevrey order $1$, but is not of any $1/q$-Gevrey order 
less than $1$.

\bigskip

Comparing Theorem \mainTh\ and Proposition \mainThPolynomial, one wonders if
Theorem \mainTh\ could be improved by showing that there is a usual
convergence and not merely in the three weak senses of 
Definition \convergenceMeaning. While we do not know if such an improvement 
is possible, the 
following result suggests that it might not be so in general, for the reason
that we will discuss after the statement. We will use the $1/q$-Borel 
transform $\calB_{1/q;1}$ introduced by B\'ezivin (1992) and Ramis (1992)
in the theory of $q$-difference equation as
$$
  \calB_{1/q;1}g(z)=\sum_{n\geq 0} q^{n\choose 2} g_n z^n\, .
$$
This transform was introduced earlier by Garsia (1981) who calls it the 
unroofing operator, but the $1/q$-Borel transform notation will be more 
convenient for our purpose.

Since $g$ solving \rightInverseDef\ also depends on $q$, we
will write $g(z;1/q)$ when we need to emphasize this dependence.
Recall that $g$ is related to $\phi$ via \rightInverseDef\ and \phiDef.

\Proposition{\label{cvToRenewal}
  Consider the sequence $(\tau_n)$ defined by $\tau_0=1$ and 
  $\tau_n=\sum_{1\leq i\leq n} \phi_i \tau_{n-i}$, whose generating function is
  $\tau(z)=1/\bigl(1-\phi(z)\bigr)$.  Then, under assumptions (i) and (ii)
  of Theorem \mainTh,

  \medskip

  \noindent
  (i) $\lim_{q\to 0} q^{n+1\choose 2} g_{n+1}=\tau_n$ for any nonnegative $n$;

  \medskip

  \noindent
  (ii) $\lim_{q\to 0} \calB_{1/q;1} g(z;1/q)/z=1/\bigl(1-\phi(z)\bigr)$.
}

\bigskip

An interesting feature of this result is that it shows a connection between
the $q$-Lagrange inversion and the probabilistic theory of renewal sequences.
Indeed, following Feller (1968; chapter XIII), we can interpret $\tau_n$
as the probability that a recurent event occurs at the $n$-th trial while
$\phi_n$ is the probability that it occurs for the first time at the $n$-th
trial.  Assertion (ii) of Proposition \cvToRenewal\ informs us that via the
$1/q$-Borel transform, the $q$-Lagrange inversion provides an interpolation
between the usual Lagrange inversion ($q=1$) and the renewal theory (limit for
$q=0+$).

Erd\"os, Pollard and Feller (1949) showed that $(\tau_n)$ converges
to $1/\sum_{n\geq 1} n\phi_n$, a quantity which may be $0$. When it does vanish,
Garsia and Lamperti (1962) showed that under condition (iii) of Theorem 
\mainTh\ and when $\rho$ is less than $1$, then
$$
  \liminf_{n\to\infty} f(1-1/n)\phi_n={\sin \pi \rho\over \pi} \, .
$$
When $\rho$ is strictly between $1/2$ and $1$, this liminf is in fact a limit. 
However when $\rho$ is at most $1/2$, they show that there is no limit in 
general. In light of Proposition \cvToRenewal.i, it is then uncertain if
the convergence in Theorem \mainTh\ could be strengthened into a usual
limit for all values of $\rho$, or only for $\rho$ strictly between $1/2$ and
$1$, or for $\rho$ in some other interval depending on $q$.

\bigskip

A variation on the proof of Theorem \mainTh\ gives the following, which deals
with an oscillating sequence $\bigl([z^n]f(z)\bigr)$ and yields an 
oscillating sequence $(g_n)$ as well.

\Theorem{\label{cor}
  Assume that the assumptions of Theorem \mainTh\ or Proposition 
  \mainThPolynomial\ hold with $\phi(-z)$ 
  substituted for $\phi(z)$. Then their conclusions hold with 
  the same
  substitution and $\bigl((-1)^n g_n\bigr)$ substituted for $(g_n)$.
}

\bigskip

Theorems \mainTh\ and \cor\ are obtained using a Tauberian theorem. Using
singularity analysis as in Flajolet and Sedgewick (2009) allows one to prove 
analogous results under slightly different assumptions that require 
analyticity of $\phi$ in sectors near $\zeta$ ---~see the discussion on 
singularity analysis versus Tauberian theorem and Darboux's method in 
Flajolet and Sedgewick, 2009, \S VI.11. However, it does not make the proof
any easier or more complicated.

\bigskip

\section{The \poorBold{$q$}-Catalan numbers} When $f(z)=z(1-z)$, that
is, when $\phi(z)=z$, equation \rightInverseDef\ becomes
$$
  g(z)-g(z)g(z/q)=z \, .
  \eqno{\equa{qCatalanGF}}
$$
The coefficients $\calC_n=q^{n+1\choose 2}g_{n+1}$ are the $q$-Catalan numbers
and were introduced by Carlitz and Riordan (1964) and Carlitz (1972).

The assumptions of Proposition \mainThPolynomial\ are satisfied, since
$\phi_n=\One\{\, n=1\,\}$, and the only root of $\phi(\zeta)=1$ is $\zeta=1$.
Therefore, Proposition \mainThPolynomial\ yields
$$
  \limn \calC_n 
  = \limn q^{n\choose 2} g_n
  = 1/(q;q)_\infty \, .
  \eqno{\equa{qCatalan}}
$$
This result was obtained by F\"urlinger, and Hofbauer (1985) using a
combinatorial argument. Further results on $q$-Catalan numbers are in Mazza 
and Piau (2002).

The core 
argument of the proof of Theorem \mainTh\ yields an interesting formula 
for some of the coefficients of the $q$-Catalan numbers.

\Proposition{\label{qCatalanCase}
  Let $\calC_n(q)$ be the $n$-th $q$-Catalan number. Let $j$ be a positive
  integer. If $n$ is at least $j$, then 
  $[q^j]\calC_n(q) = [q^j]{1/(q;q)_\infty}$.
}

\bigskip

\Proof As described in F\"urlinger and Hofbauer (1985) and as can be derived
from \qCatalanGF, the $q$-Catalan numbers obey the recursion
$$
  \calC_{n+1} 
  = \sum_{0\leq m\leq n} \calC_{n-m} \calC_m q^{m(n-m+1)} \, ,
  \eqno{\equa{qCatalanCaseA}}
$$
with $\calC_0=1$.
Consider the exponent of $q$ in the right hand side of this identity. The 
function $m\in\{\, 0,1,\ldots,n\,\}\mapsto m(n-m+1)$ is $0$ at $m=0$; otherwise
it is at least $n$. Therefore, if $m$ does not vanish, $\calC_{n-m}\calC_m
q^{m(n-m+1)}$ viewed as a polynomial in $q$ is in the ideal generated 
by $q^n$. Consequently, applying
$[q^j]$ on both sides of \qCatalanCaseA\ we obtain for $n>j$,
$$
  [q^j]\calC_{n+1}=[q^j]\calC_n \,.
  \eqno{\equa{qCatalanCaseB}}
$$
Consequently, $[q^j]\calC_n$ does not depend on $n$ provided $n$ is at 
least $j$. But Theorem \mainTh\ implies that $([q^j]\calC_n)_{n\geq 0}$ 
converges to $[q^j]1/(q;q)_\infty$ when $0<q<1$. This 
implies Proposition \qCatalanCase.\hfill\qed

\bigskip

Note that F\"urlinger and Hofbauer (1985) provide a combinatorial 
interpretation of $[q^j]\calC_n(q)$ for any $j$ and $n$, from which
they deduce that $\limn \calC_n(q)=1/(q;q)_\infty$ in the formal
topology (see Flajolet and Sedgewick, 2009, \S A.5) and this follows from
Proposition \qCatalanCase\ as well.

\bigskip


\section{Examples}
The purpose of this section is to give a couple of examples which illustrate
the results of section \fixedref{2}. Throughout this section, we assume that
$0<q<1$.

\medskip

\noindent{\it Example 1.} A simple fractional extension of the $q$-Catalan 
numbers is obtained when $f(z)=z(1-z)^\rho$ for some positive $\rho$ less 
than $1$. In this case
$$
  \phi(z)
  =1-(1-z)^\rho
  = \sum_{n\geq 1} {\rho (1-\rho)\cdots (n-1-\rho)\over n!} z^n
$$
and $\zeta=1$. We see that for any nonnegative integer $j$,
$$\eqalign{
  {\phi_{n+j}\over \phi_n}
  &{}={(n+j-1-\rho)\cdots (n-\rho)\over (n+j)\cdots (n+1)} \cr
  &{}=\Bigl( 1-{\rho+1\over n+j}\Bigr) \cdots \Bigl(1-{\rho+1\over n+1}\Bigr)
  \, .\cr
}
$$
In particular, $\limn \phi_{n+j}/\phi_n=1$ and assumption (iii) of Theorem
\mainTh\ holds. To check assumption (iv),
$$
  {1\over 1-\phi(z)} 
  = (1-z)^{-\rho}
  = 1+\sum_{n\geq 1} {\rho(\rho+1)\cdots (\rho+n-1)\over n!} z^n \, .
$$
Thus,
$$
  {[z^{n+1}]\over [z^n]}{1\over (1-\phi)}
  = {\rho+n\over 1+n} \, .
$$
Assumption (iv) of Theorem \mainTh\ then follows from a result of
Bojanic and Seneta (1973) (see Bingham, Goldie and Teugels, 1989; 
Theorem 1.9.8) We then
obtain, in the sense of the convergences in limit infimum, except on a
set of integers of density $0$, and coefficientwise,
$$
  q^{n\choose 2} n^{1-\rho} g_n
  \to {1\over\Gamma(\rho) \prod_{j\geq 1} (1-q^j)^\rho} 
  = {1\over\Gamma(\rho) (q;q)_\infty^\rho} \, .
$$
In order to have an exact analogue of \qCatalan, we would need to have a
usual convergence, but we do not know how to prove it.

\medskip

\noindent {\it Example 2.} Let $\lambda$ be a positive real number and 
consider the function $\phi(z)=e^\lambda (e^z-1)$. Thus,
$$
  f(z)=z-e^\lambda \sum_{k\geq 1} {z^{k+1}\over k!} \, .
$$
The equation $\phi(z)=1$ has a unique solution $\zeta=\log (1+e^{-\lambda})$.
Since $\phi_n=e^\lambda/n!$ for any $n\geq 1$, we see that for any 
positive integer $j$,
$$
  \limn \phi_{n+j}/\phi_n=0 \, ,
$$
and assumption (iii) of Theorem \mainTh\ is satisfied. 

To check assumption (iv) of Theorem \mainTh, we have
$$\eqalign{
  1-\phi(\zeta z)
  &{}=(1+e^\lambda)\bigl( 1-(1+e^{-\lambda})^{z-1}\bigr) \cr
  &{}\sim (1+e^\lambda) (1-z) \log (1+e^{-\lambda})
  }
$$
as $z$ tends to $1$. Setting $\epsilon=(1-z)\log (1+e^{-\lambda})$, we obtain
$$\displaylines{\qquad
  {1\over 1-\phi(\zeta z)} 
  - {1\over (1+e^\lambda) (1-z)\log(1+e^{-\lambda})}
  \hfill\cr\hfill
    = {1\over 1+e^\lambda} 
    { (e^{-\epsilon}-1+\epsilon)/\epsilon^2 \over (1-e^{-\epsilon})/\epsilon } 
    \, .
  \qquad\cr}
$$
The function
$$
  { (e^{-\epsilon}-1+\epsilon)/\epsilon^2 \over (1-e^{-\epsilon})/\epsilon }
  = {\sum_{k\geq 2} (-1)^k \epsilon^{k-2}/k!
     \over -\sum_{k\geq 1} (-1)^k\epsilon^{k-1}/k!}
$$
is continuous on the closed disk $|z|\leq 1$ and continuously differentiable
on $|z|=1$. Therefore, using Darboux's method (Flajolet and Sedgewick, 2009,
Theorem VI.14)
$$\eqalign{
  [z^n] {1\over 1-\phi(\zeta z)}
  &{}=[z^n] {1\over (1+e^\lambda)(1-z)\log(1+e^{-\lambda})} + o(1) \cr
  &{}={1\over (1+e^\lambda)\log(1+e^{-\lambda})} +o(1) \cr
  }
$$
as $n$ tends to infinity. Consequently, the sequence 
$[z^n]1/\bigl(1-\phi(\zeta z)\bigr)$, $n\geq 1$, is regularly varying of
index $\rho-1$ with $\rho=1$.  Since $\phi'(\zeta)=1+e^\lambda$, Theorem
\mainTh\ yields
$$
  q^{n\choose 2} \bigl(\log(1+e^{-\lambda})\bigr)^n (1+e^\lambda)g_n
  \to {1\over \prod_{j\geq 1} 
              \Bigl( 1-e^\lambda \bigl( (1+e^{-\lambda})^{q^j}-1\bigr)\bigr)}
$$
in the sense of convergence in limit infimum, except on a set of integers of
density $0$ and coefficientwise.

\bigskip

%

\section{A result when \poorBold{$q>1$}{\relax}}%
As before we write $f(z)=z\bigl(1-\phi(z)\bigr)$. In the range $q>1$,
following B\'ezivin (1992) and Ramis (1992), recall that $(\phi_n)$ is 
in a $q$-Gevrey class of order $s$, if there exists a constant $A$ such that
$$
  \phi_n=O( A^n q^{s{n\choose 2}})
$$
as $n$ tends to infinity. We will assume slightly more than the belonging
to the $q$-Gevrey class of order $1$, namely,
$$
  \sum_{n\geq 1} \phi_n z^n q^{-{\ss n-1\choose\ss 2}}<\infty 
  \quad\hbox{for any positive $z$.}
  \eqno{\equa{phiGevrey}}
$$
We then have the following result, which does not assume that $(\phi_n)$ is
a nonnegative sequence as we did in the previous sections.

\Theorem{\label{qGOne} 
  Assume that $q>1$ and \phiGevrey\ holds. Then $g$
  has a positive radius of convergence $\eta$ which satistifies
  $$
    \sum_{k\geq 1} \phi_k g\Bigl({\eta\over q}\Bigr)\cdots 
    g\Bigl({\eta\over q^k}\Bigr) = 1 \, .
  $$
  Furthermore, if one defines $C$ by
  $$
    {1\over C}=\sum_{k\geq 1} \phi_k g\Bigl({\eta\over q}\Bigr)\cdots 
    g\Bigl({\eta\over q^k}\Bigr) \sum_{1\leq i\leq k} {1\over q^i} 
    {g'\over g}\Bigl({\eta\over q^i}\Bigr) \, ,
  $$
  then
  $$
    g_n\sim C/\eta^n
  $$
  as $n$ tends to infinity.
}

\bigskip

Under the assumptions of Theorem \qGOne, the sequence $(g_n)$ has a remarkably
simple asymptotic behavior, though the constant is not explicit and may
be a complicated function of $f$ and $q$. In the special case of 
the $q$-Catalan numbers, that is when $\phi(z)=z$, the result
is due to Mazza and Piau (2002).

\bigskip

%
\section{Connection with linear \poorBold{$q$}-difference equations}
In this section only, we will consider first some solutions 
of \rightInverseDef\ which do not vanish at $0$.
In equation \rightInverseDef, we may make a change of function as in
Gessel (1980) or Brak and Prellberg (1995), setting for some nonzero real 
number $\eta$,
$$
  g(z)=\eta {h(z/q)\over h(z)} \, ;
$$
in particular, $g_0\not=0$. This simplifies \rightInverseDef\ to a linear
$1/q$-difference equation
$$
  \sum_{n\geq 1} f_n \eta^n h(z/q^n)=zh(z) \, .
  \eqno{\equa{qLinearA}}
$$
The general theory of linear $q$-difference equations, when only a finite 
number of coefficients do not vanish and are power series in $z$ goes back to 
Carmichael (1912), Birkhoff (1913), 
Adams (1929, 1931), and Trjitzinsky (1938) and has been revived recently by
B\'ezivin (1992), Ramis (1992), Sauloy (2000, 2003) and Zhang (1999, 2002)
among others (see Ramis, Sauloy and Zhang, 2013, for references). The modern
theory classifies these equations and studies the divergent series that are
formal solutions, parallelling the theory of ordinary differential equations
and their formal power series solutions. While this theory is 
relevant to \qLinearA,
it appears to provide only upper bound on the coefficients of the power
series $h$ when $f$ is a polynomial, and those bounds do not seem to translate
into a precise estimate on the coefficients of $g$. However, \qLinearA\
is still interesting in connection to \rightInverseDef, since it allows
us to provide explicit formal solutions to \rightInverseDef. For this
purpose, we need the following definition.

\Definition{\label{qExtremalZero}
  A zero $\kappa$ of $f$ is $q$-extremal if $f(\kappa/q^n)\not=0$ for 
  any $n\geq 1$.
}

\bigskip

We then have the following set of formal solutions to \rightInverseDef,
and it is an open problem to determine if other solutions exist.

\Proposition{\label{formalSolutions}
  For each $q$-extremal zero $\kappa$ of $f$, \rightInverseDef\ has a 
  formal solution
  $$
    g_\kappa(z)
    =\kappa {\ds 1+\sum_{n\geq 1} 
     {  \ds z^n\over{\ds q^n}\prod_{\ss 1\leq i\leq n}\ds f(\kappa/q^i)}
      \over
        \ds 1+\sum_{n\geq 1} 
        {\ds z^n\over\prod_{\ss 1\leq i\leq n}\ds f(\kappa/q^i)} }
  \, . 
  $$
}

\Proof 
Let $\kappa$ be a $q$-extremal zero of $f$. We set $\eta=\kappa$ 
in \qLinearA\ and apply $[z^k]$ to both sides of \qLinearA\ to obtain
$$
  \cases{ \ds\sum_{n\geq 1} f_n {\kappa^n\over q^{nk}} h_k =  h_{k-1}
          & if $k\geq 1$,\cr\noalign{\vskip 3pt}
          \ds\sum_{n\geq 1} f_n \kappa^n h_0=0\,,&\cr}
$$
that is
$$
  \cases{ f(\kappa/q^k) h_k=h_{k-1} & if $k\geq 1$
          \cr\noalign{\vskip 3pt}
          f(\kappa)\tilde h_0 = 0\,. &\cr}
  \eqno{\equa{formalSolutionsA}}
$$
Since $\kappa$ is $q$-extremal, $f(\kappa/q^k)$ does not vanish for any positive
$k$ and, by induction, \formalSolutionsA\ yields
$$
  h_k = {h_0\over \prod_{1\leq i\leq k} f(\kappa/q^i)} \, .
$$
Thus,
$$
  h(z)
  =h_0 \sum_{n\geq 0} {z^n\over \prod_{1\leq i\leq n} f(\kappa/q^i)}
  \eqno{\equa{formalSolutionB}}
$$
and $g$ follows.\hfill\qed

\bigskip

It is easy to give sufficient conditions for $g_\kappa$ to be an actual
solution. To fix the ideas, let us consider $q$ in $(0,1)$. 
If $\liminf_{t\to\sign(\kappa)\infty}|f(t)|>0$, we see that the function $h$ 
in \formalSolutionB\ has a positive radius of convergence and $g_\eta$
is a well defined function. Conversely, if $\lim_{t\to\sign(\kappa)\infty} f(t)=0$,
then $h$ is a divergent series and the solution $g_\kappa$ is only a 
formal one.\note{\vskip -35pt study the confluence problem; add a section on thermodynamic
of the $q$-Lagrange inversion?}

\medskip

One may wonder if one could use a similar change of function, 
namely, $g(z)=zh(z/q)/h(z)$, in order to recover the solution of
\rightInverseDef\ with vanishing constant coefficient. Such a change
of function transforms \rightInverseDef\ into
$$
  \sum_{n\geq 1} f_n {z^n\over q^{n\choose 2}} h(z/q^n)
  = zh(z) \, .
$$
Since $h$ is determined up to a multiplicative constant, we take $h_0=1$.
Applying $[z^k]$ to both sides of this identity, we obtain
$$
  \sum_{1\leq n\leq k} {f_n \over q^{{n\choose 2}+n(k-n)}} h_{k-n}
  = h_{k-1} \, ,
$$
that is, setting $m=k-n$,
$$
  \sum_{0\leq m\leq k-2} {f_{k-m} \over q^{(k-m)(k+m-1)/2}} h_m
  = h_{k-1}\Bigl( 1-{f_1\over q^{k-1}}\Bigr) \, .
$$
In the setting of Theorem \mainTh, we have $f(z)=z\bigl( 1-\phi(z)\bigr)$,
so that $f_1=1$ and $f_n=-\phi_{n-1}$. Substituting $k+1$
for $k$, we obtain
$$
  \Bigl( {1\over q^k}-1\Bigr) h_k
  = \sum_{0\leq m\leq k-1} {\phi_{k-m}\over q^{(k-m+1)(k+m)/2}} h_m 
  \, .
  \eqno{\equa{divergentSolution}}
$$

In order to analyze the asymptotic behavior of the sequence $(h_n)$, we now
distinguish according to the position of $q$ with respect to $1$.

\medskip

\noindent{\it Case $0<q<1$.} Recall that we took $h_0=1$. 
Since $1/q^k>1$ for any $k\geq 1$, 
identity \divergentSolution\ shows that if the $\phi_i$ are nonnegative, so 
are the $h_k$. Then, let $N$ be such that $\phi_N$ is positive. Isolating 
the term for which $m=k-N$ in \divergentSolution, we obtain
$$
  \Bigl({1\over q^k}-1\Bigr) h_k
  \geq {\phi_N \over q^{(N+1)(2k-N)/2}} h_{k-N} \, .
$$
Since $1/q^k\geq (1/q^k)-1$, this yields
$$
  h_k\geq {h_{k-N}\over q^{kN}} q^{N(N+1)/2} \phi_N \, .
$$
By induction, since $h_0=1$, then at least some
$h_m$ is positive, and then $h_{m+kN}$ grows at least like $q^{-N^2k^2/2(1+o(1))}$
as $k$ tends to infinity. This precludes $h(z)$ to be a convergent series,
which was to be expected since the solution of \rightInverseDef\ is a 
divergent series.

\medskip

\noindent{\it Case $q>1$.} Since $(k-m+1)(k+m)\geq 2k$ in the 
range $0\leq m\leq k-1$, identity \divergentSolution\ yields
$$
  \Bigl( 1-{1\over q^k}\Bigr) |h_k|
  \leq {1\over q^k} \sum_{0\leq m\leq k-1}\phi_{k-m} |h_m| \, .
  \eqno{\equa{convergentSolution}}
$$
We will now show that this inequality implies that $h$ has infinite radius of
convergence, as soon as the $q$-Borel 
transform $\calB_{q;1}\phi(z)=\sum_{j\geq 0} q^{-j^2/2}\phi_jz^j$ has an infinite 
radius of convergence, that is, as soon as \phiGevrey\ holds. In particular, 
if $\phi(z)$ has a positive radius of convergence, then the radius of 
convergence of $h$ is infinite.

Set $H_k=\max_{0\leq i\leq k} |h_i|$. If $0\leq m\leq k-1$, then
$(k-m+1)(k+m)\geq (k-m)^2$. Thus, identity \divergentSolution\ yields
$$
  \Bigl(1-{1\over q}\Bigr) |h_k|
  \leq \sum_{0\leq m\leq k-1} {\phi_{k-m}\over q^{(k-m)^2/2}} H_{k-1} \, .
$$
Therefore,
$$
  H_k\leq {q\over q-1} \calB_{q;1}\phi(1) H_{k-1} \, .
$$
By induction, this shows that there exists a positive constant $c$ such 
that $H_k\leq c^k$. Consequently, $\tilde h(z)=\sum_{n\geq 0} |h_n|z^n$ has 
a positive radius of convergence. 

In the range $0\leq m\leq k-1$, we also have $(k-m+1)(k+m)\geq (k-m)^2+2m$. 
Therefore, \divergentSolution\ implies
$$
  \Bigl(1-{1\over q}\Bigr) |h_k|
  \leq \sum_{0\leq m\leq k-1} {\phi_{k-m}\over q^{(k-m)^2/2}} {|h_m|\over q^m}
  \, . 
$$
Considering a positive real number $x$ and multiplying both sides of this
inequality by $x^k$, summing over $k$, setting $j=k-m-1$, and using
that $\phi_0=0$, we obtain
$$\eqalignno{
  \Bigl(1-{1\over q}\Bigr) \tilde h(x)
  &{}\leq \sum_{\ss m\geq 0\atop\ss j\geq 0} {\phi_{j+1}\over q^{(j+1)^2/2}} x^{j+1}
   {|h_m|\over q^m} x^m \cr
  &{}=\tilde h(x/q)\calB_{q;1} \phi(x) \, . 
  &\equa{convergentSolutionA}\cr
  }
$$
By assumption, the $q$-Borel transform $\calB_{q;1} \phi$ 
has infinite
radius of convergence. Since the radius of convergence of $\tilde h$ is 
positive, inequality \convergentSolutionA\ shows that this radius is infinite.
The radius of convergence of $h$ coincides with that of $\tilde h$, and so is
also infinite.

\bigskip

%

\section{Proof of Theorem \mainTh}
Define the sequence $(t_n)_{n\in\NN}$ by
$$
  \sum_{n\geq 0} t_n q^{-{n+1\choose 2}} f(z)\cdots f(q^nz)=z \, .
  \eqno{\equa{tnDef}}
$$
Applying $U_{f,q}^{-1}$, that is, $U_{g,1/q}$, on both sides of this identity 
and switching the left and right hand sides,
$$
  g(z)=\sum_{n\geq 0} t_n q^{-{n+1\choose 2}} z^{n+1} \, .
$$
Hence, $g_0=0$ and 
$$
  g_{n+1}=q^{-{n+1\choose 2}} t_n
  \eqno{\equa{gntn}}
$$
for any nonnegative integer $n$.

The proof consists in studying the sequence $(t_n)$. Because the strategy
used is likely to be useful in studying similar combinatorial sequences
related to $q$-series, and because the length of the proof may hide the
simplicity of the approach in the technical details, we now outline
the proof. We will first show that we can assume that $\zeta=1$, which
we do from now on in this outline.

It is convenient to define
$$
  \Phi_n=\Gamma(\rho)\Bigl(1-\phi\Bigl(1-{1\over n}\Bigr)\Bigr)
$$
and
$$
  L(q)={1\over \prod_{j\geq 1} \bigl(1-\phi(q^j)\bigr)} \, .
$$

Since $g_{n+1}=q^{-{n+1\choose 2}} t_n$ and $(\Phi_n)$ is regularly varying, 
it is easy to see that Theorem \mainTh\ is equivalent to the following.

\Theorem{\label{mainThtn}
  Let $(t_n)$ be as in \tnDef. Suppose that assumptions (i)---(iv) of
  Theorem \mainTh\ hold and that $\zeta=1$. Then,

  \smallskip
  \noindent
  (i) $\liminfn n\Phi_n t_n = L(q)$;

  \smallskip
  \noindent
  (ii) the sequence $(n\Phi_n t_n)_{n\geq 0}$ converges to $L(q)$ except
  possibly on a set of integers of density $0$;

  \smallskip
  \noindent
  (iii) for any $j$,  $\limn n\Phi_n [q^j]t_n=[q^j]L(q)$.
}

\bigskip

To continue our outline, the first step to prove Theorem \mainThtn\ is 
to use a form of generating function for the sequence
$(t_n)$. In our case, it turns out that relation \tnDef\ and representation
\phiDef\ yield the identity
$$
  \bigl(1-\phi(z)\bigr) \sum_{n\geq 0} t_n z^n \prod_{1\leq j\leq n} 
  \bigl(1-\phi(q^jz)\bigr) = 1 \, .
$$
In this identity, for $n$ large, the products are about
$\prod_{j\geq 1} \bigl(1-\phi(q^jz)\bigr)$, which, when $z$ is $1$ is $1/L(q)$. 
This suggests that
$$
  \lim_{z\to 1}\bigl(1-\phi(z)\bigr) \sum_{n\geq 0} t_n z^n 
  = L(q) \, .
$$
From such a limit, a Tauberian theorem gives immediately that
$$
  \limn \rho\Phi_n \sum_{0\leq i<n} t_i 
  = L(q) \, .
  \eqno{\equa{heuristicA}}
$$
This is like a Ces\`aro limit of the sequence $(t_n)$. Since there does not seem
to be a way to show that the sequence $(t_n)$ is ultimately monotone, we need
to resort to further arguments to go from this Ces\`aro type limit to an
actual limit. However, the standard argument based on monotonicity suggests
that we should have $t_n\sim C/\rho n\Phi_n$ as $n$ tends to infinity. In 
this case we should expect \heuristicA\ to yield
$$
  \limn \Phi_n \sum_{0\leq i<n} {C\over i\Phi_i}  = L(q) \, .
$$
Since $(\Phi_n)$ is regularly varying of negative index $-\rho$, Karamata's 
theorem implies
$$
  \limn \Phi_n \sum_{0\leq i<n} {1\over i\Phi_i}= {1\over \rho} \, .
$$
Thus, $C=\rho L(q)$ and we would have $\limn n\Phi_nt_n=L(q)$. 

To justify this heuristic, we need a more direct description of the sequence
$(t_n)$. This is provided by a recursion which at first looks hopelessly 
complicated. To write it, for any tuple of nonnegative integers 
$(n_1,\ldots,n_i)$ with
sum $n_1+\cdots+n_i=n$, we define
$$
  L_i(n_1,\ldots,n_i)
  ={n\choose 2}-\sum_{1\leq j\leq i} {n_j\choose 2}
  + \sum_{1\leq j\leq i} (j-1)n_j \, .
  \eqno{\equa{LiDef}}
$$
Then we will see that $(t_n)$ obeys the recursion
$$\displaylines{\qquad
    t_n
    =\One\{\, n=0\,\} +\sum_{i\geq 1} \phi_i \sum_{n_1,\ldots,n_{i+1}}
    t_{n_1}\cdots t_{n_{i+1}} q^{L_{i+1}(n_1,\ldots,n_{i+1})}
  \hfill\cr\hfill
    \One\{\, n_1+\cdots+n_{i+1}=n-i\,\} \, . 
  \qquad\equa{Segner}\cr}
$$

How may one study such a recursion once one knows that $(t_n)$ converges in some
form of Ces\`aro sense? As indicated in Barbe and MacCormick (2012), recursion
\Segner\ is a generalization of \qCatalanCaseA. In particular \Segner\
implies that each $t_n$ is a polynomial in $q$; thus as we did in the proof 
of Proposition \qCatalanCase, we may consider $[q^j]t_n$. We will see that
$L_{i+1}(n_1,\ldots,n_{i+1})=0$ if and only if $n_2=\cdots=n_{i+1}=0$. Thus,
defining
$$
  \epsilon_n = t_n-\sum_{1\leq i\leq n} \phi_i t_{n-i} \, ,
  \eqno{\equa{epsilonDef}}
$$
recursion \Segner\ implies that $\epsilon_n$ is a sum of polynomials in the 
ideals generated by the $q^{L_{i+1}(n_1,\ldots,n_{i+1})}$ with $n_1<n-i$. 
It happens that $L_{i+1}(n_1,\ldots,n_{i+1})$ is fairly large for most 
tuples $(n_1,\ldots,n_{i+1})$. This suggests that if we consider 
$[q^j]\epsilon_n$, very few terms will remain ---~in fact about $j^j$
terms, a large number in terms of $j$, but a small one in terms of $n$. 
This is the key to estimate the magnitude of 
$[q^j]\epsilon_n$. Then, we can express $t_n$ in terms of $\epsilon_n$ and show
that for any integer $j$, the sequence $(n\Phi_n[q^j]t_n)_{n\geq 0}$ converges;
we can then identify this limit as $[q^j]L(q)$.
By summing over $j$, these limits provide the inequality 
$\liminfn n\Phi_n t_n\geq L(q)$. This and the Ces\`aro type limit yield 
Theorem \mainThtn. We now turn to the actual proof, which we will divide into
several steps for the sake of clarity.

\step{Step 1. Reduction to \poorBold{$\zeta$}=\hskip1pt 1}%
Since the proof is slightly easier to write 
when $\zeta$ is $1$, we show in this preliminary step how to reduce the 
general case to that where $\zeta$ is $1$.

\Lemma{\label{reductionLemma}
  If $g(z)$ is the right inverse of $f(z)$, then $g(\lambda z)/\lambda$ is
  the right inverse of $f(\lambda z)/\lambda$.
}

\bigskip

\Proof 
Let $\tilde g(z)$ be the right inverse of $f(\lambda z)/\lambda$. 
Since $[z^n]\bigl(f(\lambda z)/\lambda\bigr)$ is $\lambda^{n-1} f_n$, we deduce 
from \rightInverseDef\ that
$$
  \sum_{n\geq 1} f_n\lambda^{n-1} \tilde g(z)\tilde g(z/q)\cdots 
  \tilde g(z/q^{n-1}) = z \, ,
$$
that is
$$
  \sum_{n\geq 1} f_n \bigl(\lambda \tilde g(z)\bigr) 
  \bigl(\lambda \tilde g(z/q)\bigr) \cdots 
  \bigl( \lambda\tilde g(z/q^{n-1})\bigr) =\lambda z \, .
$$
It then follows from the uniqueness of the right inverse and the definition
\rightInverseDef\ that $\lambda \tilde g(z) =g(\lambda z)$.\hfill\qed

\bigskip

We now explain why Theorem \mainTh\ is equivalent to Theorem \mainThtn. Start
with $f$ and $\phi$ as in Theorem \mainTh, that is with $\zeta$ which may not
be $1$. Set $\tilde f(z)=f(\zeta z)/\zeta$. This change of function ensures
that $[z^1]\tilde f=[z^1]f$. The smallest nonnegative
zero of $\tilde f$ is $\tilde\zeta=1$. Writing 
$\tilde f(z)=z\bigl(1-\tilde \phi(z)\bigr)$, we see that 
$\tilde\phi(z)=\phi(\zeta z)$. Thus,
$$
  1-\tilde \phi(1-1/n))=1-\phi(\zeta(1-1/n)\bigr)
$$
and
$$
  \prod_{j\geq 1} \bigl( 1-\tilde\phi(q^j)\bigr)
  = \prod_{j\geq 1} \bigl( 1-\phi(\zeta q^j)\bigr) \, .
$$
Since $g(z)=\zeta\tilde g(z/\zeta)$, we have $g_n=\zeta^{1-n}\tilde g_n$.
The equivalence between Theoren \mainTh\ and \mainThtn\ then follows.

\bigskip

To conclude this first step, we can assume that $\zeta=1$, which we do 
from now on.

\bigskip

\step{Step 2. Ces\`aro type limit for \poorBold{$(t_n)$}} 
Our goal in this subsection is to prove the following result.

\Proposition{\label{CesaroLimittn}
  $\limn \rho \Phi_n\sum_{0\leq i<n} t_i = L(q)$.
  \hfill{\rm\equa{tCesaro}}
}

\bigskip

\Proof It follows from \tnDef\ of this paper, Theorem 5.3 
and identity (5.4) in Barbe and McCormick (2012) that $(t_n)$
may be interpreted as dual coefficients associated with the Catalan power
series $P(z,t)=t-t\phi(tz)$. Thus, setting $\phi_i=[z^i]\phi(z)$, 
Theorem 3.3 in Barbe and McCormick (2012) and its translation as 
identities (5.18) and (5.19) in that paper imply that $(t_n)$ satisfies 
the recursion \Segner\ --- alternatively, this can be proved directly 
from \rightInverseDef\ as indicated in the remark following this proof.
Since $(\phi_i)_{i\geq 1}$ is a nonnegative sequence under assumption (i) 
and $t_0=1$, an induction shows that all the $t_i$'s are nonnegative when
$q$ is nonnegative.

Define the power series
$$
  H(z)=\sum_{n\geq 0} t_n z^n \prod_{1\leq j\leq n} \bigl( 1-\phi(q^jz)\bigr)
  \, .
  \eqno{\equa{HDef}}
$$
Given \phiDef, \tnDef\ is
$$
  \bigl( 1-\phi(z)\bigr) H(z)=1 \, . 
  \eqno{\equa{basic}}
$$

At this point, $H$ is a formal power series. To make it a convergent one,
define for any integer
$n$,
$$
  A_n=t_n\prod_{1\leq j\leq n} \bigl( 1-\phi(q^j)\bigr) \, ,
$$
with $A_0=t_0=1$, and set
$$
  A(x)=\sum_{n\geq 0} A_n x^n \, .
$$
Let $x$ be a positive real number less than $\zeta=1$. 
Using assumption (i) of Theorem \mainTh, $\phi$ is 
nondecreasing and $1-\phi(q^{j+1}x)\geq 1-\phi(q^{j+1})$. Thus, using 
that $\phi$ is nonnegative under assumption (i),
$$
  H(x)
  \geq \sum_{n\geq 0} t_n x^n 
    \prod_{1\leq j\leq n}\bigl( 1-\phi(q^j)\bigr)
  = A(x) \, .
$$
This inequality, combined with \basic, shows that $A(x)$ is a convergent series
when $x$ is in $[\,0,1)$, and
$$
  \limsup_{x\uparrow 1}\bigl( 1-\phi(x)\bigr) A(x)
  \leq 1 \, .
  \eqno{\equa{eqLimSup}}
$$

Next, we also have
$$\eqalign{
  H(x)
  &{}=\prod_{j\geq 1} \bigl( 1-\phi(q^j x)\bigr)
    \sum_{n\geq 0} {t_n x^n\over \prod_{j> n} \bigl( 1-\phi(q^j x)\bigr)} \cr
  &{}\leq \prod_{j\geq 1} \bigl( 1-\phi(q^j x)\bigr)
    \sum_{n\geq 0} {t_n x^n\over \prod_{j> n} 
    \bigl( 1-\phi(q^j)\bigr)} \cr
  &{}\leq { \prod_{j\geq 1} \bigl( 1-\phi(q^j x)\bigr)\over
            \prod_{j\geq 1} \bigl( 1-\phi(q^j)\bigr) } A(x) \, . \cr}
$$
Therefore, \basic\ yields
$$
  \liminf_{x\uparrow 1} \bigl( 1-\phi(x)\bigr) A(x)
  \geq 1 \, .
$$
Combined with \eqLimSup, this implies
$$
  \lim_{x\uparrow 1} \bigl( 1-\phi(x)\bigr) A(x)
  = 1 \, .
  \eqno{\equa{AEquiv}}
$$
In other words, $A(x)\sim 1/\bigl(1-\phi(x)\bigr)$ as $x$ tends to 1. As
indicated after \mainThA, assumption (iv) implies that 
the function $t\mapsto A(1-1/t)$ is regularly varying at
infinity with index $\rho$. Karamata's Tauberian theorem for power series
(see Bingham, Goldie and Teugels, 1989, Corollary 1.7.3) and \AEquiv\ imply 
that $\sum_{0\leq i<n} A_i \sim 1/\rho\Phi_n$ as $n$ tends to infinity. Since
$A_i\sim t_i /L(q)$ as $i$ tends to infinity, we obtain \tCesaro.

\Remark In the proof of Proposition \CesaroLimittn\ we showed \Segner\ 
by interpreting the sequence $(t_n)$ as dual coefficients of some $q$-Catalan
basis. This can also be shown directly from \rightInverseDef\ as follows.
Since $g_0=0$, we set $g(z)=zh(z)$, and since $f_0=1$, we 
rewrite \rightInverseDef\ as
$$
  zh(z)
  =z-\sum_{n\geq 2} {f_n\over q^{n\choose 2}} z^n h(z)h(z/q)\cdots h(z/q^{n-1}
  \, .
$$
Since $f_n=-\phi_{n-1}$ for any $n\geq 2$,
$$
  h(z)=1+\sum_{i\geq 1} {\phi_i\over q^{i+1\choose 2}} z^i 
  h(z)h(z/q)\cdots h(z/q^i) \, .
$$
Applying $[z^n]$ we obtain
$$\displaylines{\quad
  h_n=\One\{\, n=0\,\}+\sum_{i\geq 0} {\phi_i\over q^{i+1\choose 2}}
  \sum_{n_{i+1},n_i,\ldots,n_1} q^{-(0n_{i+1}+1n_i+\cdots+in_1)}
  \hfill\cr\hfill
  h_{n_{i+1}}h_{n_i}\cdots h_{n_1} \One\{\, n_1+\cdots +n_{i+1}=n-i\,\}\, .
  \cr}
$$
Since $t_n=q^{n+1\choose 2}g_{n+1}=q^{n+1\choose 2} h_n$, we obtain
$$\displaylines{\quad
  t_n=\One\{\, n=0\,\} +\sum_{i\geq 0} \phi_i \sum_{n_1,\ldots,n_{i+1}}
  q^{L_{i+1}(n_1,\ldots,n_{i+1})}
  \hfill\cr\hfill
  t_{n_1}\cdots t_{n_{i+1}}
  \One\{\, n_1+\cdots+n_{i+1}=n-i\,\}
  \cr}
$$
with now
$$\displaylines{\quad
  L_{i+1}(n_1,\ldots,n_{i+1})
  = -\sum_{1\leq j\leq i+1} (i+1-j)n_j
  \hfill\cr\hfill
  {}-\sum_{1\leq j\leq i+1}{n_j+1\choose 2}
  + {n+1\choose 2}-{i+1\choose 2} \, .
  \qquad\equa{LiAltDef}\cr}
$$
and $n_1+\cdots+n_{i+1}+i+1=n+1$. Therefore, if $n_1+\cdots+n_i=n$,
\LiAltDef\ yields
$$
  L_i(n_1,\ldots,n_i)=-\sum_{1\leq j\leq i} (i-j)n_j -\sum_{1\leq j\leq i} {n_j+1\choose 2} +{n+i\choose 2} -{i\choose 2} \, .
  \eqno{\equa{LiAltDefA}}
$$
Note that $\ds{n_j+1\choose 2}={n_j\choose 2}+n_j$, so that \LiAltDefA\ is
$$
  -(i-1)n +\sum_{1\leq j\leq i}(j-1)n_j -\sum_{1\leq j\leq i} {n_j\choose 2}-n
  {}+{n+i\choose 2}-{i\choose 2} \, .
$$
After an elementary calculation, this coincides with \LiDef.

\step{Step 3. Some results on \poorBold{$L_i$}} As indicated previously, if the 
sequence $(t_i)$ were monotone, Theorems \mainThtn\ and \mainTh\ would follow
readily from Proposition \CesaroLimittn. 
However, it does not seem that there is a simple argument
to show that $(t_n)$ is monotone, even ultimately. To proceed further,
we need to use whatever extra information we have on $(t_n)$, namely the 
explicit description given by recursion \Segner. In this description, the key
observation to make is that if $L_{i+1}(n_1,\ldots,n_{i+1})$ is at least
a given integer $j$ it is conceivable that the right hand side of \Segner\
is much simpler when applying $[q^j]$. It is of course a leap of faith,
but the example of the $q$-Catalan numbers suggests that it is not hopeless.
This raises immediately the question of studying the quantity
$L_{i+1}(n_1,\ldots,n_{i+1})$ and understand reasonably well its order
of magnitude in relation to the tuple $(n_1,\ldots,n_{i+1})$. 

Note that if $(n_1,\ldots,n_i)$ is a tuple with $n=n_1+\cdots+n_i$, 
then, writing explicitly the binomial coefficients involved in $L_i$, we 
obtain
$$
  L_i(n_1,\ldots,n_i)
  = {n(n-2)\over 2} -{1\over 2} \sum_{1\leq j\leq i} n_j^2
  +\sum_{1\leq j\leq i} j n_j \, .
  \eqno{\equa{LiAlt}}
$$
To study $L_i$ we will need the following definition.

\Definition{\label{basicTransform}
  Given two positive integers $j$ and $k$ and $j<k$,

  \medskip
  \noindent 
  (i) the raising operator $R_{j,k}$ acts on tuples of length $i\geq k$ by 
  increasing the $j$-th component of the tuple by $1$ and decreasing
  the $k$-the component by $1$;

  \medskip
  \noindent
  (ii) the transposition $\tau_{j,k}$ acts on a tuple of length $i\geq k$ by
  permuting $n_j$ and $n_k$.
}

\bigskip

For instance, we have
$$
  R_{2,5}(1,2,3,4,5,6,7)=(1,3,3,4,4,6,7) \, 
$$
and
$$
  \tau_{2,5}(1,2,3,4,5,6,7)=(1,5,3,4,2,6,7) \, .
$$
Raising operators are defined in Macdonald (1995, \S I.1).

Note that the raising operators do not change the sum of the elements of
tuples. The following result implies that raising operators decrease the 
value of $L_i$ on nonincreasing
tuples. The second assertion of the following lemma is essentially in Hardy,
Littlewood and Polya (1952, \S 10.2) but we reproduce it here for convenience.

\Lemma{\label{basicTransformL}
  Let $j<k\leq i$ be some positive integers. If $n_k\geq 1$ then
  $$
    L_i(n_1,\ldots,n_i)=L_i\circ R_{j,k}(n_1,\ldots,n_i)+n_j-n_k-j+k+1
  $$
  and
  $$
    L_i(n_1,\ldots,n_i)=L_i\circ \tau_{j,k}(n_1,\ldots,n_i)+(k-j)(n_k-n_j) \, .
  $$
}

\noindent
The sole purpose of the condition $n_k\geq 1$ is to ensure 
that $R_{j,k}(n_1,\ldots, n_i)$ is a tuple of nonnegative integers. The result 
holds under more general conditions.

\bigskip

\Proof Since the raising operators do not change the sum of the elements
of a tuple, \LiAlt\ yields
$$\displaylines{\quad
  L_i(n_1,\ldots,n_i)-L_i\circ R_{j,k}(n_1,\ldots,n_i)
  =-{1\over 2} (n_j^2+n_k^2)+jn_j+kn_k 
  \hfill\cr\hfill
  {}+ {1\over 2}\bigl( (n_j+1)^2+(n_k-1)^2\bigr) -j(n_j+1)-k(n_k-1) \, ,\cr
}
$$
which, after some simplification is the first assertion.

Since the transpositions change neither the sum of the elements of a tuple
nor the sum of their square, \LiAlt\ yields
$$
  L_i(n_1,\ldots,n_i)-L_i\circ \tau_{j,k}(n_1,\ldots,n_i)
  = jn_j+kn_k -jn_k-kn_j \, ,
$$
which is the second assertion.\hfill\qed

\bigskip

In the first assertion of Lemma \basicTransformL, we see that the term
$n_j-n_k+(k-j)+1$ is positive whenever $j<k$ and $n_j\geq n_k$. Hence, as
announced, applying a raising operator to a nonincreasing tuple decreases the 
corresponding value of $L_i$. This implies the following result which tells
us when $L_i$ is small and whose third assertion informs us on the sparsity of
tuples for which $L_i$ exceeds a specified value.

\Lemma{\label{LBound}
  (i) $\min_{n_1+\cdots+n_i=n} L_i(n_1,\ldots,n_i)=0$, the minimum being achieved
  at the unique tuple $(n,0,\ldots,0)$. 

  \medskip

  \noindent (ii) $\min_{\matrix{\ss n_1+\cdots+n_i=n
                                \hfill\cr\noalign{\vskip -4pt}
                                \ss n_1<n\hfill\cr}} L_i(n_1,\ldots,n_i)
  = n$, the minimum being achieved at the unique tuple $(n-1,1,0,\ldots,0)$.

  \medskip

  \noindent (iii) If $n_1+\cdots+n_i=n$ and $k$ is the greatest index for which
  $n_k$ is positive, then $L_i(n_1,\ldots,n_i)\geq n+k-2$.
}

\bigskip

\Proof (i) Consider a tuple $(n_1,\ldots,n_i)$ with $n_1+\cdots+n_i=n$.
By ordering its elements in nonincreaing order we minimize the sum $\sum jn_j$
involved in \LiAlt. Thus assume from now on that $n_1\geq n_2\geq\ldots 
\geq n_i$. If some
$n_j$ is positive, we can decrease $L_i$ by applying the operator
$R_{1,j}$.

\noindent (ii) The same argument applies for this second assertion, 
except that since 
$n_1$ is at most $n-1$, we must keep some $n_j$ with $j\geq 2$ positive. The
best we can do after ordering the tuple in nonincreasing order is by 
applying some raising operators $R_{2,j}$ as many times as possible 
and then the raising operator $R_{1,2}$ until we reduce $n_2$ to $1$.
We then obtain the tuple $(n-1,1,0,\ldots,0)$ for which $L_i$ has, from
its definition, the value ${n\choose 2}-{n-1\choose 2}+1$, which is $n$.

\noindent (iii) Consider a tuple for which $k$ is the greatest index for
which $n_k$ does not vanish. We can decrease the value of $L_i$ on this tuple
using transpositions to order, in nonincreasing order, only the positive 
entries of the tuples; for instance, we would 
reorder $(6,0,0,0,5,7,0,2,0,0)$ as $(7,0,0,0,6,5,0,2,0,0)$. Then, we can
reduce further the value of $L_i$ by applying some
raising operator to bring the tuple to the form $(n-1,0,\ldots,0,1,0,\ldots,0)$
where the $1$ entry is the $k$-th. Since
$$
  L_i(n-1,0,\ldots,0,1,0,\ldots, 0)
  ={n\choose 2}-{n-1\choose 2}+k-1
  = n+k-2 \, ,
$$
this proves (iii).\hfill\qed

\bigskip

We now make a parenthetical comment concerning this third step and the 
combinatorics of partitions and tableaux. Recall that a composition of $n$ is
a tuple $(n_1,\ldots, n_i)$ of nonnegative integers which sum to $n$. A
partition of $n$ is a composition $(n_1,\ldots, n_i)$ of $n$ such that
$n_1\geq \cdots\geq n_i$. These tuples can be extended to sequences, by
agreeing that $n_j$ is $0$ if $j>i$. We refer to Macdonald (1995) for the 
terminology that we will use.

For a partition $(n_1,\ldots,n_i)$, one defines the conjugate partition
$n'_j=\sharp\{\, i\,:\, n_i\geq j\,\}$. We then have
$$
  \sum_{i\geq 1} (i-1)n_i=\sum_{j\geq 1} {n'_j\choose 2} \, .
$$
Thus,
$$
  L_i(n_1,\ldots,n_i)
  ={n\choose 2}-\sum_{i\geq 1} {n_i\choose 2}+\sum_{j\geq 1} {n'_j\choose 2}
  \, .
  \eqno{\equa{stepThreeA}}
$$
Lemma \LBound\ asserts that $L_i(n_1,\ldots,n_i)$ is minimum when $(n_1,\ldots,
n_i)$ is $(n,0,\ldots,0)$, and we can also check that $L_i(n_1,\ldots,n_i)$
is maximal when $(n_1,\ldots,n_i)$ is $(0,0,\ldots,n)$, the maximal value being
$n(i-1)$. The difference
$$
  \Delta(n_1,\ldots,n_i)
  = \sum_{j\geq 1} {n_j\choose 2}-\sum_{i\geq 1}{n'_i\choose 2}
$$
involved in \stepThreeA\ vanishes when $n_i=n'_i$ and is a measure of the
discrepency between a partition and its conjugate, or analogously of the
asymmetry of the tableau associated to a partition. This discrepency is in
some sense maximal when the partition is $(n,0,\ldots,0)$ and \stepThreeA\
is
$$
  \Delta(n_1,\ldots,n_i)-\Delta(n,0,\ldots,0) \, .
$$
Thus, by calculating this difference,
$L_i$ measures a form of closeness to the most unsymmetric partition
$(n,0,\ldots,0)$.

\step{Step 4. Bounding \poorBold{$[\,q^j]\epsilon_n$}} The purpose of this part 
of the
proof is to obtain some reasonable bounds on $[q^j]\epsilon_n$ which we will
then use for showing that $(n\Phi_n [q^j]t_n)_{n\geq 0}$ converges.
In the course of the proof we will need the operators $P_k$ which map
a power series $g(z)=\sum_{i\geq 0} g_i z^i$ to its projection on the space
spanned by $z^0, z^1, \ldots, z^{k-1}$, that is
$$
  P_k g(z)= \sum_{0\leq i<k} g_i z^i \, ,
$$
with the convention that $P_0$ is $0$.

\Lemma{\label{epsilonBound} 
  The following hold:

  \medskip

  \noindent (i) $[q^0]\epsilon_n = \One\{\, n=0\,\}$. 

  \medskip

  \noindent (ii) For any nonnegative integer $j$ there exists a nonnegative
  real number $c_j$ such that $[q^j]\epsilon_n\sim c_j \phi_{n-j}$ as $n$ tends
  to infinity.

  \medskip

  \noindent (iii) For any nonnegative integer $j$ there exists a real
  number $\tilde c_j$ such that for any $n$,
  $$
    0\leq [q^j]\epsilon_n \leq \tilde c_j \sum_{n-j\leq i\leq n-1} \phi_i \, .
  $$
}

Our proof does not excludes that $c_j$ in the second assertion of the lemma 
may vanish. If this is the case, the statement should be read as
$[q^j]\epsilon_n=o(\phi_{n-j})$.

\bigskip

\Proof In the recursion \Segner, taking into consideration the second 
indicator function, that $t_0=1$ and Lemma \LBound.i holds, if $n_1=n-i$ then
$$
  t_{n_1}\ldots t_{n_{i+1}}q^{L_{i+1}(n_1,\ldots,n_{i+1})}= t_{n-i} \, .
$$
Moreover, if $i=n$ then the only tuple $(n_1,\ldots,n_{i+1})$ such that
$n_1+\cdots+n_{i+1}=n-i$ is that for which $n_1=n-i=0$ and all the other $n_j$
vanish as well. Therefore, isolating the term for which $n_1=n-i$ in 
recursion \Segner, we obtain, given how $\epsilon_n$ is defined in \epsilonDef,
$$\displaylines{\quad
  \epsilon_n 
  = \One\{\, n=0\,\} +\sum_{1\leq i\leq n-1} \phi_i
  \sum_{\matrix{\ss n_1,\ldots,n_{i+1}\hfill
                \cr\noalign{\vskip -4pt}
                \ss n_1<n-i\hfill\cr}} t_{n_1}\cdots t_{n_{i+1}}
  q^{L_{i+1}(n_1,\ldots,n_{i+1})}
  \hfill\cr\hfill
  \One\{\, n_1+\cdots+n_{i+1}=n-i\,\} \, .
  \qquad\equa{epsilonBoundA}\cr}
$$

\noindent (i) If $n\geq 1$ and $i\leq n-1$ and $n_1+\cdots+n_{i+1}=n-i$ and
$n_1<n-i$, then Lemma \LBound.ii shows that
$$
  L_{i+1}(n_1,\ldots,n_{i+1})\geq n-i\geq 1\, .
  \eqno{\equa{epsilonBoundB}}
$$
Thus, if $n$ is at least $1$, each 
term $t_{n_1}\cdots t_{n_{i+1}}q^{L_{i+1}(n_1,\ldots,n_{i+1})}$ involved 
in \epsilonBoundA\ is in the
ideal generated by $q$. This implies $[q^0]\epsilon_n=\delta_{0,n}$.

\hfuzz=2pt
\noindent (ii) Given assumption (iii) of Theorem \mainTh, 
set $\rho_j=\limn \phi_{n+j}/\phi_n$.  Then, define
$$\displaylines{\quad
  c_j=\sum_{1\leq i\leq j} \rho_{j-i}
      \sum_{\matrix{\ss n_1,\ldots,n_{j+1}\hfill
                    \cr\noalign{\vskip -4pt}
                    \ss n_1<i\hfill\cr}}
  [q^{j-L_{j+1}(n_1,\ldots,n_{j+1})}](t_{n_1}\cdots t_{n_{j+1}})
  \hfill\cr\noalign{\vskip -10pt}\hfill
  \One\{\, n_1+\cdots+n_{j+1}=i\,\} \, .
  \quad\cr}
$$
\hfuzz=0pt
Since \epsilonBoundB\ holds, the term
$t_{n_1}\cdots t_{n_{i+1}} q^{L_{i+1}(n_1,\ldots,n_{i+1})}$ in \epsilonBoundA\
is in the ideal
generated by $q^{n-i}$. Thus, whenever $n-i>j$, 
$$
  [q^j]( t_{n_1}\cdots t_{n_{i+1}}q^{L_{i+1}(n_1,\ldots,n_{i+1})})=0\, .
$$
Therefore, applying $[q^j]$ to both sides of \epsilonBoundA, we obtain
$$\displaylines{\qquad
  [q^j]\epsilon_n
  = \sum_{n-j\leq i\leq n-1}\phi_i \sum_{\matrix{\ss n_1,\ldots,n_{i+1}\hfill
                                              \cr\noalign{\vskip -4pt}
                                              \ss n_1<n-i\hfill\cr}}
  [q^{j-L_{i+1}(n_1,\ldots,n_{i+1})}](t_{n_1}\cdots t_{n_{i+1}})
  \hfill\cr\hfill \One\{\, n_1+\cdots+n_{i+1}=n-i\,\} \, .
  \qquad\equa{epsilonBoundC}\cr
}
$$
In the range $n-i\leq j$, a tuple $(n_1,\ldots,n_{i+1})$ 
with $n_1+\cdots +n_{i+1}=n-i\leq j$ may have at most $j$ nonzero elements. 
By Lemma \LBound.iii, if the last nonzero element has index $k$ then,
for $i$ in the range of the summation in \epsilonBoundC,
$$
  L_{i+1}(n_1,\ldots,n_{i+1})\geq n-i+k-2\geq k-1 \, .
$$
Thus, if $L_{i+1}(n_1,\ldots,n_{i+1})\leq j$, then $(n_1,\ldots,n_{i+1})$
is in fact $(n_1,\ldots,n_{j+1},0,\ldots ,0)$ and
$$
  L_{i+1}(n_1,\ldots,n_{i+1}) 
  = L_{j+1}(n_1,\ldots , n_{j+1}) \, .
$$
Given \epsilonBoundC, this implies that $[q^j]\epsilon_n$ is
$$\displaylines{\quad
  \sum_{n-j\leq i\leq n-1} \phi_i 
  \sum_{\matrix{\ss n_1,\ldots,n_{j+1}\hfill
                  \cr\noalign{\vskip -4pt}
                  \ss n_1<n-i\hfill\cr}}
    [q^{j-L_{j+1}(n_1,\ldots,n_{j+1})}](t_{n_1}\cdots t_{n_{j+1}})
  \hfill\cr\noalign{\vskip -10pt}\hfill
  \One\{\, n_1+\cdots+n_{j+1}=n-i\,\}
  \qquad\cr\noalign{\vskip 15pt}\qquad\qquad
  {}=\sum_{1\leq i\leq j} \phi_{n-i} 
    \sum_{\matrix{\ss n_1,\ldots,n_{j+1}\hfill
                 \cr\noalign{\vskip -4pt}
                 \ss n_1<i\hfill\cr}}
    [q^{j-L_{j+1}(n_1,\ldots,n_{j+1})}](t_{n_1}\cdots t_{n_{j+1}})
  \hfill\cr\noalign{\vskip -10pt}\hfill
  \One\{\, n_1+\cdots+n_{j+1}=i\,\} \, . 
  \qquad\equa{epsilonBoundCa}\cr}
$$
Since the $t_n$ are polynomials in $q$ with nonnegative coefficients
and the sequence $(\phi_n)$ satisfies assumption (iii) of Theorem \mainTh, 
we obtain that $[q^j]\epsilon_n\sim \phi_{n-j} c_j$ as $n$ tends to 
infinity.

\noindent (iii) The summation over $n_1,\ldots,n_{j+1}$ and $n_1<i$ in 
\epsilonBoundCa\ is at most
$$\displaylines{\quad
  \tilde c_j=\max_{1\leq i\leq j}
  \sum_{\matrix{\ss n_1,\ldots,n_{j+1}\hfill
                 \cr\noalign{\vskip -4pt}
                 \ss n_1<i\hfill\cr}}
    [q^{j-L_{j+1}(n_1,\ldots,n_{j+1})}](t_{n_1}\cdots t_{n_{j+1}})
  \hfill\cr\noalign{\vskip -10pt}\hfill
  \One\{\, n_1+\cdots+n_{j+1}=i\,\} \, . \quad\cr}
$$
The third assertion of the lemma follows.\hfill\qed

\bigskip

It will prove convenient to record the following consequences of 
Lemma \epsilonBound.

\Lemma{\label{epsilonSumBound}
  (i) For any $1\leq j\leq n_1< n_2$,
  $$
    \sum_{n_1\leq n< n_2}[q^j]\epsilon_n 
    \leq \tilde c_j (n_2-n_1) \sum_{n_1-j\leq i<n_2} \phi_i \, .
  $$

  \medskip

  \noindent (ii) $[q^j]\epsilon_n\leq j\tilde c_j\max_{n-j\leq i\leq n-1}\phi_i$.

  \medskip

  \noindent (iii) $\sum_{i\geq 0} [q^j]\epsilon_i \leq j\tilde c_j$.
}

\bigskip

In particular, assertion (iii) of Lemma \epsilonSumBound\ shows that the 
series of nonnegative terms $\sum_{i\geq 0} [q^j]\epsilon_i$ is finite.

\bigskip

\Proof (i) Lemma \epsilonBound.iii\ ensures that the left hand side in the 
statement is at most
$$\displaylines{\quad
  \tilde c_j \sum_{n,i}\One\{\, n-j\leq i\leq n-1\,;\, n_1\leq n< n_2\,\}\phi_i
  \hfill\cr\hfill
  {}\leq \tilde c_j \sum_{n,i} 
  \One\{\, n_1-j\leq i< n_2-1 \,;\, n_1\leq n< n_2\,\}
  \phi_i \, \quad\cr}
$$
This upper bound is the right hand side of the statement.

\noindent (ii) It follows immediately from Lemma \epsilonBound.iii.

\noindent (iii)Using Lemma \epsilonBound.iii,
$$\eqalign{
  \sum_{n\geq 0} [q^j]\epsilon_n
  &{}\leq \tilde c_j\sum_{n,i}
          \One\{\, n-j\leq i\leq n-1\,;\, n\geq 0\,\} \phi_i \cr
  &{}=\tilde c_j \sum_{n,i} \One\{\, i+1\leq n\leq i+j\,;\, n\geq 0\,\} 
      \phi_i \cr
  &{}\leq \tilde c_j j\sum_i \phi_i \, . \cr}
$$
The result follows since $\sum_{i\geq 1} \phi_i=\phi(1)=1$.\hfill\qed

\bigskip

\step{Step 5. \poorBold{$(n\Phi_n [\,q^j]t_n)_{n\geq 0}$} converges} We are 
now in position
to prove that $(n\Phi_n[q^j]t_n)$ is a convergent sequence, even though
we will not be able to explicitly describe its limit in this step.

\Lemma{\label{tnCv}
  For any integer $j$,
  $$
    \limn n\Phi_n [q^j]t_n=\sum_{i\geq 0} [q^j]\epsilon_n
  $$
  and this limit is finite.
}

\bigskip

We will show later how this limit can be made explicit in terms of the original
data of the problem, namely the function $\phi$.

To prove Lemma \tnCv, we will need the following auxilliary result.
Recall that we write $[z^i]1/(1-\phi)$ for the coefficient of $z^i$ in 
the power expansion of $1/\bigl(1-\phi(z)\bigr)$.

\Lemma{\label{tnCvLemma}
  For any nonnegative integer $j$,
  $$
    \lim_{N\to\infty} \limsup_{n\to\infty} n\Phi_n \sum_{0\leq i\leq n-N}
    [z^i]{1\over 1-\phi} \phi_{n-j-i} = 0 \, .
  $$
}

\bigskip

\Proof We first prove the assertion when $j$ vanishes. 
Note that $\phi_{n-i}=[z^{n-i}]\phi$ and for any $n\geq 1$
$$\eqalign{
  \sum_{0\leq i\leq n} [z^i]{1\over 1-\phi}[z^{n-i}]\phi
  &{}= [z^n]{\phi\over 1-\phi}
   {}= [z^n]\Bigl( -1+{1\over 1-\phi}\Bigr) \cr
  &{}= [z^n]{1\over 1-\phi} \, . \cr}
$$
Therefore, for any $n\geq N$,
$$
  \sum_{0\leq i\leq n-N} [z^i]{1\over 1-\phi} [z^{n-i}]\phi
  = [z^n]{1\over 1-\phi} -\sum_{n-N<i\leq n} [z^i]{1\over 1-\phi} \phi_{n-i}
  \, . \eqno{\equa{tnCvA}}
$$
Using assumption (iv) of Theorem \mainTh, which implies \mainThA, 
uniformly in $i$ between $n-N$ and $n$,
$$
  [z^i]{1\over 1-\phi} \sim [z^n]{1\over 1-\phi} \sim {1\over n\Phi_n}
  \eqno{\equa{tnCvB}}
$$
as $n$ tends to infinity. Thus, \tnCvA\ is
$$
  {1\over n\Phi_n} \Bigl( 1+o(1)-\sum_{0\leq i<N}\phi_i\bigl(1+o(1)\bigr)\Bigr)
  \, .
$$
The assertion follows since $\sum_{i\geq 0} \phi_i=\phi(1)=1$.

Assume now that $j$ is positive. Set
$$
  A_{n,j,N}=n\Phi_n \sum_{0\leq i\leq n-N} [z^i]{1\over 1-\phi} \phi_{n-j-i}
$$
Since 
$$
  A_{n,j,N}
  = {n\Phi_n\over (n-j)\Phi_{n-j}} A_{n-j,0,N-j}\, ,
$$
the result follows from the regular variation of the sequence $(n\Phi_n)$
and the result when $j$ vanishes.\hfill\qed

\bigskip

\noindent{\bf Proof of Lemma \tnCv.}
 Let $B$ be the backward shift, mapping $t_n$ to $t_{n-1}$. We agree to set
$t_n=0$ if $n$ is negative, so that $\epsilon_n=(1-\phi)(B)t_n$. We then have
$$
  t_n
  ={1\over 1-\phi}(B)\epsilon_n
  = \sum_{0\leq i\leq n} \Bigl([z^i]{1\over 1-\phi}\Bigr) \epsilon_{n-i} \, .
  \eqno{\equa{tnCvAa}}
$$
Applying $[q^j]$ to both sides of this identity and multiplying by
$n\Phi_n$,
$$
  n\Phi_n [q^j]t_n 
  = n\Phi_n \sum_{0\leq i\leq n} [z^i]{1\over 1-\phi} [q^j]\epsilon_{n-i} \, .
$$
We will prove that the right hand side converges and this will require
us to split the sum into two parts.

Since $1/(1-\phi)=\sum_{k\geq 0} \phi^k$, and the $\phi_i$ are nonnegative,
all the coefficients $[z^i]\bigl(1/(1-\phi)\bigr)$ are also nonnegative. 

Let $\eta$ be a positive real number less than $1$. Using
Lemma \epsilonBound.ii, there exists an integer $N$ such that for any
integer $n$ at least $N$, the coefficient $[q^j]\epsilon_n$ is between
$(1-\eta)c_j\phi_{n-j}$ and $(1+\eta)c_j\phi_{n-j}$ when $c_j$ is positive. 
Thus, whenever $n$ is at least $N$,
$$\displaylines{\qquad
  n\Phi_n \sum_{0\leq i<n-N} [z^i] {1\over 1-\phi} [q^j]\epsilon_{n-i} 
  \hfill\cr\hfill
  {}\leq (1+\eta) c_j n\Phi_n \sum_{0\leq i<n-N}[z^i]{1\over 1-\phi} \phi_{n-j-i}
  \, .
  \qquad\equa{tnCvAb}\cr}
$$
Using Lemma \tnCvLemma, this upper bound can be made arbitrarily small by
taking $N$ large and $n$ large enough. Substituting $-\eta$ for $\eta$ gives
a matching lower bound. The case $c_j=0$ is handled analogously.

Next, using \tnCvB, and the nonnegativity of $[q^j]\epsilon_i$, we also obtain
$$
  n\Phi_n \sum_{n-N\leq i\leq n} [z^i] {1\over 1-\phi} [q^j]\epsilon_{n-i}
  \sim \sum_{0\leq i\leq N} [q^j]\epsilon_i \, ,
  \eqno{\equa{tnCvC}}
$$
as $n$ tends to infinity. Combining \tnCvAb\ and \tnCvC, we obtain that
$$
  \limn n\Phi_n [q^j]t_n
  = \sum_{i\geq 0} [q^j]\epsilon_i
$$
and the limiting sum is finite as indicated in 
Lemma \epsilonSumBound.\hfill\qed

\bigskip

\step{Step 6. Ces\`aro limit for (\poorBold{$[\,q^j]t_n$})}
The drawback of Lemma \tnCv\ is that it does not tell us what the limit
of $(n\Phi_n[q^j]t_n)_{n\geq 0}$ is. To identify this limit,
our next result shows that we can apply $[q^j]$ to both sides of \tCesaro\ 
and permute it with the limit in the left hand side.

\Lemma{\label{tCoeffCesaro}
  For any nonnegative integer $j$,
  $$
    \limn \rho\Phi_n\sum_{0\leq i<n} [q^j]t_i 
    = [q^j]L(q) \, . 
  $$
}

Before proving Lemma \tCoeffCesaro, we need an elementary result. Recall
that whenever $k$ is an integer, $P_k$
is the operator that maps a power series $\sum_{i\geq 0} g_i x^i$ to its 
projection on the space spanned by $x^0,\ldots, x^{k-1}$
with the convention that $P_0$ is $0$.

\Lemma{\label{truncateSeries}
  Let $g(z)=\sum_{i\geq 0} g_i z^i$ be a series with nonnegative coefficients. 
  Then, for any nonnegative integer $k$ and any real numbers $0\leq x\leq y$,
  $$
    0
    \leq g(x)-P_k g(x)
    \leq \Bigl({x\over y}\Bigr)^k g(y) \, .
  $$
}

\Proof It follows from
$$
  0
  \leq \sum_{i\geq k} g_i x^i
  \leq \sum_{i\geq k} \Bigl({x\over y}\Bigr)^i  g_i y^i
$$
and that $(x/y)^i\leq (x/y)^k$ whenever $i$ is at least $k$.\hfill\qed

\bigskip

\noindent{\bf Proof of Lemma \tCoeffCesaro.} 
We write $\overline t_n$ for $\rho\Phi_n\sum_{0\leq i<n} t_i$.
Recursion \Segner\ shows that $t_n=t_n(q)$ is a polynomial in $q$ with 
nonnegative coefficients. In particular if $q$ is at most $r$ then,
using Proposition \CesaroLimittn,
$$
  q^k [q^k]\overline t_n 
  \leq \overline t_n(q)
  \leq\overline t_n(r)
  \leq L(r)+o(1)
$$
where the $o(1)$ term is as $n$ tends to infinity and uniform in the index $k$.
Consequently, for any $n$ large enough, 
$([q^k]\overline t_n)_{k\geq 0}$ is in $\prod_{k\geq 0} [\,0,2L(r)/q^k\,]$. By
Tychonoff's theorem, the sequence of sequences $n\mapsto 
([q^k]\overline t_n)_{k\geq 0}$ belongs to a compact set for the product 
topology. It has a cluster point. Let $(a_k)_{k\geq 0}$ be such a cluster
point and let $n_j \mapsto ([q^k]\overline t_{n_j})_{k\geq 0}$ be a subsequence
which converges to that cluster point as $j$ tends to infinity. 
Set $a(x)=\sum_{k\geq 0} a_k x^k$. Since the $a_k$'s are necessarily 
nonnegative, the function $a$ is well defined on $[\,0,\infty)$, 
though possibly infinite. However, since $a_k$ is in $[\,0,2L(r)/q^k\,]$, 
the radius of convergence of $a$ is at least $q$ as long as $L(r)$ is finite.
By lemma \truncateSeries,
$$
  0\leq t_i(q)-P_kt_i(q)
  \leq \Bigl({q\over r}\Bigr)^k t_i(r) \, .
$$
Thus, summing over $i$ from $0$ to $n-1$ and multiplying by $\rho\Phi_n$,
$$
  0\leq \overline t_n(q)-P_k\overline t_n(q)
  \leq \Bigl({q\over r}\Bigr)^k \overline t_n(r) \,,
$$
and, taking limit along the subsequence $(n_j)$ and using 
Proposition \CesaroLimittn, 
$$
  0\leq L(q)-P_k a(q)\leq \Bigl({q\over r}\Bigr)^k L(r) \, .
$$ 
It follows that $\lim_{k\to\infty} P_k a(q)=L(q)$ for any $q$ less than the 
radius of convergence of $L$. Since $P_k a(q)$ is the partial 
sum $\sum_{0\leq n< k} a_nq^n$, this means that $a$ and $L$ coincide. 
Consequently, the
sequence $n\mapsto ([q^k]\overline t_n)_{k\geq 0}$ has a unique cluster point,
which is $([q^k]L)_{k\geq 0}$, and it converges to this cluster point,
which is what Lemma \tCoeffCesaro\ asserts.\hfill\qed

\bigskip

\step{Step 7. Identifying the limit of (\poorBold{$n\Phi_n[\,q^j]t_n$})}
Following the principle that a sequence which converges converges
to its Ces\`aro limit, we combine Lemmas \tnCv\ and \tCoeffCesaro\ to
obtain the following result.

\Lemma{\label{tnLimit}
  $\limn n\Phi_n [q^j]t_n= [q^j]L$.
}

\bigskip

\Proof Using Lemma \tnCv, the sequence $(n\Phi_n [q^j]t_n)_{n\geq 0}$ 
converges to a limit $M$ which depends on $j$. We need to show that $M$ 
is $[q^j]L$. Given Lemma
\tCoeffCesaro, we consider $\Phi_n\sum_{0\leq i<n} [q^j]t_i$. Let $N$ be
a positive integer. Since $(\Phi_n)$ converges to $0$,
$$
  \limN\limn\Phi_n\sum_{0\leq i\leq N} [q^j]t_i=0 \, .
$$
Next, let $\epsilon$ be a positive real number. We have
$$
  \Phi_n\sum_{N\leq i< n} [q^j]t_i
  = \Phi_n \sum_{N\leq i<n} {1\over i\Phi_i} 
    i\Phi_i [q^j]t_i \, .
$$
If $N$ is large enough, this sum is at most
$$
  (1+\epsilon) \Phi_n M \sum_{N\leq i<n} {1\over i\Phi_i} \,,
$$
and at least the same bound with $-\epsilon$ substituted for $\epsilon$. 
Using Karamata's theorem, this sum is equivalent to
$$
  (1+\epsilon) \Phi_n M {1\over \rho\Phi_n} \, .
$$
Since $\epsilon$ is arbitrary this shows that
$$
  \limn \Phi_n \sum_{0\leq i<n} [q^j]t_n
  =M/\rho \, .
$$
Lemma \tCoeffCesaro\ then yields $[q^j]L=M$.\hfill\qed

\step{Concluding the proof of Theorem \mainThtn}
We first prove part of assertion (i) of Theorem \mainThtn.
Let $j$ be a nonnegative integer. Since $t_n\geq P_jt_n$, Lemma
\tnLimit\ yields
$$
  \liminfn n \Phi_n t_n
  \geq \limn \Phi_n P_j t_n
  = P_j L \, .
$$
Since $j$ is arbitrary, and $\lim_{j\to\infty}P_j L=L$, we proved that
$$
  \liminfn n\Phi_n t_n\geq L \, .
  \eqno{\equa{conclusionA}}
$$
The proof of assertion (ii) of Theorem \mainThtn\ will imply that the 
inequality in \conclusionA\ is in fact an equality.

\bigskip

\noindent{\it Proof assertion (ii) of Theorem \mainThtn.}
Our next result, which is in a more general setting than that used so far, 
shows that Proposition \CesaroLimittn\ and \conclusionA\ imply the 
convergence of $(n\Phi_n t_n)$ on a set
of integers of density $1$. The convergence except on a set of integers 
of density $0$ in Theorem \mainTh\ follows readily
since $(t_n)$ and $(g_n)$ are related through \gntn.
Moreover, with \conclusionA\ this implies 
that $\liminf_{n\to\infty} n\Phi_n t_n= L$.

\Lemma{\label{denseCv}
  Let $(t_n)$ be any nonnegative sequence such that
  
  \medskip

  \noindent
  (i) $\limn \rho\Phi_n \sum_{0\leq i<n} t_i =L$, and

  \medskip

  \noindent
  (ii) $\liminfn n\Phi_n t_n\geq L$.

  \medskip

  \noindent
  Then $(n\Phi_nt_n)$ converges to $L$ on a set of integers
  of density $1$.
}

\bigskip

\Proof Let $\epsilon$ and $\delta$ be some positive real numbers. For any
$n$ large enough,
$$\eqalignno{
  L+\epsilon
  &{}\geq \rho\Phi_n\sum_{0\leq i<n} t_i \cr
  &{}\geq \rho\Phi_n \sum_{0\leq i<N} t_i
   + \rho\Phi_n \sum_{N\leq i<n} {1\over i\Phi_i} i\Phi_i t_i 
   \One\{\, i\Phi_i t_i \leq L +\delta\,\}\cr
  &\qquad\quad {}+(L+\delta)\rho\Phi_n \sum_{N\leq i<n} {1\over i\Phi_i}
    \One\{\, i\Phi_i t_i>L +\delta\,\} \, .
  &\equa{denseCvA}\cr}
$$
Using assumption (ii), whenever $N$ is large enough,
$i\Phi_i t_i\geq L-\epsilon$ for any $i\geq N$. Then, dropping the first
sum in the right hand side of \denseCvA, \denseCvA\ is at least
$$\displaylines{\qquad
  (L-\epsilon)\rho\Phi_n \sum_{N\leq i<n} {1\over i\Phi_i}
  \One\{\, i\Phi_it_i\leq L+\delta\,\}
  \hfill\cr\hfill
    {}+(L+\delta)\rho\Phi_n \sum_{N\leq i<n} {1\over i\Phi_i}
    \One\{\, i\Phi_i t_i >L+\delta\,\} \,.\qquad\cr}
$$
Collecting the terms involving $L$, this is at least
$$\displaylines{\qquad
    L\rho\Phi_n \sum_{N\leq i<n} {1\over i\Phi_i}
    -\epsilon\rho\Phi_n \sum_{N\leq i<n} {1\over i\Phi_i}
  \hfill\cr\hfill
    {}+{\delta\rho\over n} \sum_{N\leq i<n} {n\Phi_n\over i\Phi_i} 
    \One\{\, i\Phi_i t_i>L +\delta\,\} \, .
    \qquad\equa{denseCvB}\cr}
$$
Since $(n\Phi_n)$ is a regularly varying sequence of index $1-\rho$ less than
$1$, Karamata's theorem implies
$$
  \sum_{N\leq i<n}{1\over i\Phi_i}\sim {1\over \rho\Phi_n}
$$
as $n$ tends to infinity. Moreover, using the uniform convergence theorem
for regularly varying sequence, for any $n$ large enough,
$$\displaylines{\qquad
  \sum_{N\leq i\leq n} {n\Phi_n\over i\Phi_i} \One\{\, i\Phi_it_i>L+\delta\,\}
  \hfill\cr\hfill
  {}\geq (1-\epsilon)\sum_{\epsilon n\leq i\leq n} 
  \Bigl({n\over i}\Bigr)^{1-\rho} \One\{\, i\Phi_i t_i>L+\delta\,\} \, .
  \qquad\cr}
$$
Since $1-\rho$ is nonnegative, this is at least
$$\displaylines{\qquad
  (1-\epsilon)\sum_{\epsilon n\leq i\leq n} \One\{\, i\Phi_i t_i>L+\delta\,\}
  \hfill\cr\hfill
  {}\geq (1-\epsilon)\sum_{0\leq i\leq n} \One\{\, i\Phi_i t_i>L+\delta\,\}
  - (1-\epsilon)\epsilon n \, .\qquad\cr
  }
$$
Thus \denseCvB\ is at least
$$
  L-\epsilon +{\delta\rho(1-\epsilon)\over n}
  \sum_{0\leq i<n} \One\{\, i\Phi_i t_i>L+\delta\,\} -(1-\epsilon)\epsilon
  +o(1)
$$
as $n$ tends to infinity. Since $\epsilon$ is arbitrary, this and \denseCvA\
yield
$$
  L
  \geq L+\delta\rho\limsupn {1\over n} \sum_{0\leq i<n} 
    \One\{\, i\Phi_i t_i>L+\delta\,\} \, ,
$$
which, after removing $L$ from both sides, is the result\hfill\qed

\bigskip

\Remark Given that the last summation in \denseCvB\ is a weighted density,
one may wonder if we could prove a result stronger than stated. This is
not the case, for, using that $(n\Phi_n)$ is regularly varying with
nonnegative index $1-\rho$,
\hfuzz=1pt
$$\displaylines{
    \Phi_n \sum_{0\leq i<n}{1\over i\Phi_i} 
    \One\{\, i\Phi_i t_i\geq L+\delta\,\}
  \hfill\cr\hfill
  \eqalign{
    {}\leq{}&{1\over n} \sum_{0\leq i<\epsilon n} {n\Phi_n\over i\Phi_i} 
             +{1\over n} \sum_{\epsilon n\leq i<n} 
             \Bigl({n\over i}\Bigr)^{1-\rho} 
             \One\{\, i\Phi_i t_i>L+\delta\,\} \bigl( 1+o(1)\bigr) \cr
    {}\leq{}&\int_0^\epsilon {\d u\over u^{1-\rho}}\bigl(1+o(1)\bigr)
             + {1\over \epsilon^{1-\rho}} {1\over n} \sum_{0\leq i<n}
             \One\{\, i\Phi_i t_i >L+\delta\,\} \bigl(1+o(1)\bigr)\cr}
   \cr}
$$
\hfuzz=0pt
as $n$ tends to infinity. Since $\epsilon$ is arbitrary, this shows that
$$
  \limn {1\over n} \sum_{0\leq i<n} \One\{\, i\Phi_i t_i>L+\delta\,\}
  = 0
  \eqno{\equa{densityA}}
$$
implies
$$
  \limn \Phi_n\sum_{0\leq i<n} {1\over i\Phi_i} 
  \One\{\, i\Phi_it_i\geq L+\delta\,\} =0\, .
  \eqno{\equa{densityB}}
$$
Thus, under assumption (iv) of Theorem \mainTh, \densityA\ and \densityB\ are
equivalent.
\bigskip

\section{Proof of Proposition \mainThPolynomial} As in the proof of Theorem
\mainTh, we assume that $\zeta$ is $1$. Assumption (ii) of 
Proposition \mainThPolynomial\ implies that $\phi$ is differentiable
at $1$, and therefore, $\rho=1$ and
$$
  \Phi_n\sim \phi'(1)/n \, .
  \eqno{\equa{polyAa}}
$$

The next lemma shows that $(\epsilon_n)$ defined in \epsilonDef\
converges to $0$ at an exponential rate.

\Lemma{\label{epsilonPolynomial}
  Under the assumption of Proposition \mainThPolynomial, for any $q$ 
  in $[\,0,1)$, there exists $r$ in $[\,0,1)$ such that $\epsilon_n\leq r^n$
  for any $n$ large enough.
}

\bigskip

\Proof As in \epsilonBoundA, we have for any $n$ at least $1$,
$$\displaylines{\quad
  \epsilon_n
  =\sum_{1\leq i\leq n-1} \phi_i q^{n-i}
  \sum_{\matrix{\ss n_1,\ldots,n_{i+1}\hfill\cr\noalign{\vskip -4pt}
                \ss n_1<n-i\hfill\cr}}
  t_{n_1}\cdots t_{n_{i+1}} q^{L_{i+1}(n_1,\ldots,n_{i+1})-(n-i)}
  \hfill\cr\noalign{\vskip -7pt}\hfill
  \One\{\, n_1+\cdots+n_{i+1}=n-i\,\} \, .
  \qquad\equa{polyB}\cr}
$$ 
Lemma \LBound\ shows that in this sum, the exponent $L_{i+1}(n_1,\ldots,n_{i+1})
-(n-i)$ is nonnegative. Then, since $q$ is in $[\,0,1)$, Proposition 
\CesaroLimittn\ and \polyAa\ yield for any $n$,
$$
  \max_{0\leq i\leq n} t_i 
  \leq \sum_{0\leq i\leq n} t_i
  \leq c(n+1) \, .
$$
Thus, \polyB\ implies
$$\qquad\eqalign{
  \epsilon_n
  &{}\leq \sum_{1\leq i\leq n} \phi_i q^{n-i} \bigl( c(n+1)\bigr)^{i+1}\cr
  \noalign{\vskip -3pt}
  &{}\hskip 72pt
    \sharp\{\, (n_1,\ldots,n_{i+1})\,:\, n_1+\cdots+n_{i+1}=n-i\,\} \cr
  \noalign{\vskip 7pt}
  &{}\leq \sum_{1\leq i\leq n} \phi_i q^{n-i} \bigl( c(n+1)\bigr)^{i+1} n^i 
    \, . \cr}
$$
Thus, for $n\geq 1$,
$$
  \epsilon_n\leq n^2\max_{1\leq i\leq n} \phi_i (cn^2)^i q^{n-i} \, .
$$
Let $r$ be in $(q,1)$. Using assumption (ii) of Proposition \mainThPolynomial, 
let $M$ be a number such that 
$$
  \limsup_{n\to\infty} {\log\phi_n\over n\log n} <-M<-2 \, .
$$
Let $i_0$ be such that $\log\phi_i\leq -Mi\log i$ whenever $i$ is 
at least $i_0$. If $i\leq i_0$, 
$$
  \phi_i (cn^2)^i q^{n-i}
  \leq (cn^2)^{i_0} q^{-i_0} \max_{1\leq j\leq i_0} \phi_j q^n \, .
$$
If $n$ is large enough, this upper bound is less than $r^n$.

If $i\geq i_0$, then 
$$
  \phi_i (cn^2)^i q^{n-i} 
  \leq \exp\bigl( -Mi\log i +i\log(cn^2)+(n-i)\log q\bigr) \, .
$$
The function
$$
  k(x)=-Mx\log x +x\log(cn^2)+(n-x)\log q
$$
has a maximum at $x^*=(cn^2/q)^{1/M}/e$ and $k(x^*)\sim n\log q$ as $n$ tends 
to infinity. Thus,
$$
  \epsilon_n\leq n^2\exp\Bigl( n(\log q)\bigl(1+o(1)\bigr)\Bigr)
$$
as $n$ tends to infinity and the result follows.\hfill\qed

\bigskip

Using Theorem IV.10 or VI.12 in Flajolet and Sedgewick (2009),
$$
  [z^n] {1\over 1-\phi}\sim {1\over \theta (1)}
$$
as $n$ tends to infinity. To conclude the proof of 
Proposition \mainThPolynomial, recall representation \tnCvAa. We substitute
$\epsilon_n$ for $[q^j]\epsilon_n$ in the proof of Lemma \tnCv, which then
implies that
$$
  \limn n\Phi_nt_n=\sum_{i\geq 0}\epsilon_i
  \eqno{\equa{polyC}}
$$
Lemma \epsilonPolynomial\ ensures that the limiting series converges.
Assumption (ii) implies that $\phi$ is differentiable at $1$, so that
$\rho$ is $1$ and $n\Phi_n\sim \phi'(1)$ as $n$ tends to infinty. Thus,
\polyC\ states that $(t_n)$ converges. Its  limit coincides with its
Ces\`aro limit, and Proposition \CesaroLimittn\ asserts that the limit is
$L(q)/\phi'(1)$. Proposition \mainThPolynomial\ then follows from relation
\gntn\ between $g_n$ and $t_n$.\hfill\qed

\bigskip

%

\section{Proof of Proposition \cvToRenewal}
Write $\epsilon_n(q)$ for what we wrote $\epsilon_n$ in \epsilonBoundA. Given
\epsilonBoundB\ and \epsilonBoundA,
$$
  \lim_{q\to 0} \epsilon_n(q)
  =\One\{\, n=0\,\} \, .
$$
Thus, \epsilonDef\ yields, with $t_n(q)$ for what we wrote $t_n$, for 
any $n\geq 1$,
$$
  \lim_{q\to 0} t_n(q)-\sum_{1\leq i\leq n} \phi_i t_{n-i}(q)=0 \, .
$$
Since each $t_n(q)$ is a polynomial in $q$, $\lim_{q\to 0} t_n(q)=t_n(0)$
exists. Hence, with $(\tau_n)$ as in Proposition \cvToRenewal, $t_n(0)=\tau_n$
for any $n\geq 1$. We also have $t_0(0)=1/f_1=1=\tau_0$.
Assertion (i) of Proposition \cvToRenewal\ then follows from \gntn.

To prove assertion (ii), note that under assumption (ii) of Theorem \mainTh, 
$\phi$ has a positive radius of convergence. For $z$ real, nonnegative and less
than $\zeta$, \HDef\ and \basic
$$
  \prod_{j\geq 1} \bigl( 1-\phi(q^jz)\bigr)\sum_{n\geq 0} t_n z^n
  \leq H(z)={1\over 1-\phi(z)}
  \leq \sum_{n\geq 0} t_n z^n \, .
  \eqno{\equa{tBound}}
$$
Since $\phi(0)=0$,
$$
  \lim_{q\to 0} \prod_{j\geq 1} \bigl( 1-\phi(q^jz)\bigr) = 1 \, .
$$
This and \tBound\ imply
$$
  \limsup_{q\to 0} \sum_{n\geq 0} t_n z^n
  \leq {1\over 1-\phi(z)}
  \leq \liminf_{q\to 0} \sum_{n\geq 0} t_n z^n \, .
$$
Given \gntn, the result follows.

\bigskip

%

\section{Proof of Theorem \cor}%
Consider still \tnDef. We still have relation \gntn, namely, 
$g_{n+1}=q^{-{n+1\choose 2}}t_n$. Set $s_n=(-1)^n t_n$. We then 
rewrite \Segner\ as
$$\displaylines{\qquad
    (-1)^n s_n
    = \One\{\, n=0\,\} +\sum_{i\geq 1} \phi_i \sum_{n_1,\ldots,n_{i+1}}
    (-1)^{n-i} s_{n_1}\cdots s_{n_{i+1}} 
    \hfill\cr\noalign{\vskip 3pt}\hfill
    q^{L_{i+1}(n_1,\ldots,n_{i+1})}
    \One\{\, n_1+\cdots+n_{i+1}=n-i\,\} \, ,\cr}
$$
that is,
$$\displaylines{
  s_n
    = \One\{\, n=0\,\} +\sum_{i\geq 1} (-1)^i\phi_i \sum_{n_1,\ldots,n_{i+1}}
    s_{n_1}\cdots s_{n_{i+1}} q^{L_{i+1}(n_1,\ldots,n_{i+1})}
  \hfill\cr\hfill
    \One\{\, n_1+\cdots+n_{i+1}=n-i\,\} \ , .\qquad\cr}
$$
Set $\widetilde\phi(z)=\phi(-z)=\sum_{i\geq 1} (-1)^i \phi_i z^i$. The
function $\widetilde\phi$ fulfills the assumptions of Theorem \mainTh.
Since $\widetilde\phi(\zeta)=1$, we have
$$
  \liminfn q^{n\choose 2} \zeta^{n+\rho-1} n\Gamma(\rho)
  \bigl( 1-\tilde\phi(\zeta-1/n)\bigr) s_n 
  = {1\over \prod_{j\geq 1} \bigl( 1-\tilde\phi(\zeta q^j)\bigr)} \, ,
$$
and the other assertions similar to those in Theorem \mainTh\ follow in the
same way.

\bigskip


\section{Proof of Theorem \qGOne}
The proof has two steps. First we make a change of function, setting
$$
  g(z)=zh(z/q)/h(z) \, .
  \eqno{\equa{qGOneA}}
$$
We then prove that \rightInverseDef\ has such a solution and that the radius
of convergence of $h$ is infinite. In the second step we then show that
$(\eta-z)g(z)$ is analytic in a neighborhood of $\eta$. Then the result
follows by singularity analysis.

\noindent{\it Step 1.} With the change of function \qGOneA, 
equation \rightInverseDef\ becomse
$$
  \sum_{k\geq 1} {f_k\over q^{k\choose 2}} z^k h\Bigl({z\over q^k}\Bigr)
  = zh(z) \, .
  \eqno{\equa{qGOneB}}
$$
Since $h$ is defined through \qGOneA\ up to a multiplicative constant and
since $h_0\not=0$ for $[z^1]g\not=0$, we may take $h_1=1$.
We divide both sides of \qGOneB\ by $z$ and apply $[z^n]$ to obtain
for $n\geq 1$
$$
  \sum_{1\leq k\leq n+1} {f_k\over q^{k\choose 2}} 
  {h_{n+1-k}\over q^{k(n+1-k)}} = h_n \, .
  \eqno{\equa{qGOneC}}
$$
Since $f_1=1$, we obtain
$$
  \sum_{2\leq k\leq n+1} {f_k h_{n+1-k}\over q^{k(2n+1-k)/2}}
  = h_n \Bigl( 1-{1\over q^n}\Bigr) \, .
$$
Set $H_n=\max_{1\leq i\leq n} |h_i|$. Identity \qGOneC\ implies
$$\eqalign{
  |h_n|\Bigl( 1-{1\over q^n}\Bigr)
  &{}\leq \sum_{2\leq k\leq n+1} {|f_k|\over q^{k(2n+1-k)/2}} H_{n-1} \cr
  &{}={H_{n-1}\over q^{2n-1}} \sum_{2\leq k\leq n+1} 
    {|f_k|\over q^{(k-2)(2n-1-k)/2}} \cr
  &{}\leq {H_{n-1}\over q^{2n-1}} \sum_{2\leq k\leq n+1} 
    {|f_k|\over q^{k-2\choose 2}} \, . \cr}
$$
This implies
$$
  |h_n|
  \leq {H_{n-1}\over q^{2n-1}} \sum_{k\geq 1} {|f_k|\over q^{k-2\choose 2}}
  \eqno{\equa{qGOneD}}
$$
and the series in this upper bound is finite under \phiGevrey. This implies
that for $n$ large enough, $|h_n|\leq H_{n-1}$, and, consequently, the sequence
$(H_n)$ is bounded. Then, \qGOneD\ implies that $(q^nh_n)$ is also bounded.
Consequently, the series $h(z)$ has a positive radius of convergence, and this
radius is at least $1/q$.

Assume now that we proved that $h$ converges in some interval $[\,0,R)$. Let
$x$ be in $[\,R,qR)$. Since $h_0=1$,
$$
  \Bigl| {f_k\over q^{k\choose 2}}x^kh\Bigl({z\over q^k}\Bigr)\Bigr|
  \sim {|f_k|\over q^{k\choose 2}} x^k h_0 
$$
as $k$ tends to infinity.
Given \phiGevrey, this implies that the left hand side of \qGOneB\ is finite.
Therefore, so is its right hand side. This proves that if $h$ converges on
$[\,0,R)$ then it converges on $[\,0,qR)$. Since $R$ may be chosen positive,
in fact, at least $1/q$, this shows that $h$ converges on the whole real line.
Its radius of convergence is infinite.

\noindent{\it Step 2.} Given \qGOneA, we see that $g$ is defined as long as
$h$ does not vanish. Consider the function
$$
  A(z)=\sum_{k\geq 1} \phi_k g(z/q)\cdots g(z/q^k) \, .
  \eqno{\equa{qGOneE}}
$$
Since $g(z/q^k)\sim g_1 z/q^k$ as $k$ tends to infinity, assumption \phiGevrey\
ensures that $A(z)$ is well defined as long as $g(z/q)$ is, that is, as long
as $z/q$ is less than the radius of convergence of $g$. Since the 
sequence $(g_n)$ is nonnegative (see the proof of the result when $q<1$ where
we had that $g_n$ was related to $t_n$ and that $t_n$ is nonnegative), the
sequence $([z^n]A)$ is also nonnegative.

Note that $\eta$, as defined in the statement of the theorem, is the smallest 
positive solution of the equation $A(\eta)=1$.
Since $f(z)=z\bigl( 1-\phi(z)\bigr)$, we can rewrite
\rightInverseDef\ as
$$
  g(z)\bigl( 1-A(z)\bigr) = z \, .
  \eqno{\equa{qGOneF}}
$$
This identity and the fact that $A$ has a positive radius of convergence
show that the radius of convergence of $g$ 
is $\eta$. The radius of convergence of $A$ is then $q\eta$.

Since the radius of convergence of $A$ is $q\eta$, which is greater 
than $\eta$, the function $A$ is analytic at $\eta$. Moreover, 
since the sequence $([z^n]A)$ is nonnegative, the derivatives
$A'(\eta)$ does not vanish. This implies that the function
$$
  {\eta-z\over A(\eta)-A(z)}
$$
is analytic. Then, \qGOneF\ yields
$$
  g(z)={1\over 1-(z/\eta)}\, {(\eta-z)\over A(\eta)-A(z)} {z\over \eta} \, .
$$
It then follows from Darboux's method or singularity analysis (see for instance
Theorem VI.12 in Flajolet, Sedgewick)  that
$$
  g_n\sim {1\over A'(\eta)\eta^n}
$$
as $n$ tends to infinity. After some elementary algebra to calculate $A'(z)$
from \qGOneE, this is the result.\hfill\qed

\bigskip


\section{Proof of Proposition \gnPolynomial}%
Since $f_1=1$, we rewrite \rightInverseDef\ as
$$
  g(z)=z-\sum_{n\geq 2} f_n g(z)g(z/q)\dots g(z/q^{n-1}) \, .
$$
Applying $[z^n]$ shows by induction on $n$ that $g_n$ is indeed a
polynomial in $1/q$. Since recursion \Segner\ shows that 
$t_n=q^{n+1\choose 2}g_{n+1}$ is a polynomial in $q$, the degree of $g_n$ as
a polynomial in $1/q$ is at most ${n\choose 2}$. This proves the first
assertion of the proposition.

Lemma \LBound\ implies that the only tuple
$(n_1,\ldots,n_{i+1})$ whose entries add to $n-i$ and for 
which $L_{i+1}(n_1,\ldots,n_{i+1})$ vanishes is $(n-i,0,\ldots,0)$. 
Therefore, \Segner\ yields
$$
  [q^0]t_n = \One\{\, n=0\,\} +\sum_{i\geq 1} \phi_i [q^0]t_{n-i} \, .
$$
This implies that the generating function $\sum_{n\geq 0} [q^0]t_nz^n=[q^0]t(z)$
satisfies
$$\eqalign{
  [q^0]t(z)
  &{}=1+[q^0]\sum_{n\geq 0} \sum_{i\geq 1} \phi_i t_{n-i}z^n \cr
  &{}=1+\phi(z)[q^0]t(z) \, .\cr
  }
$$
Therefore,  $[q^0]t(z)=1/\bigl(1-\phi(z)\bigr)$ and
$$
  [q^0]t_n = [z^n]{1\over 1-\phi(z)} \, .
$$
The result follows since $[q^0]t_n=[1/q^{n+1\choose 2}] g_{n+1}$. 

The last assertion of the proposition follows from recursion \Segner\ which
shows by induction on $n$ that $t_n$ is a polynomial in $q$ whose 
coefficients are all nonnegative whenever the $\phi_i$ are 
nonnegative.\hfill\qed

\bigskip


\section{Proof of Proposition \gDivergent}%
Let us consider first the case where $(\phi_n)$ is a nonnegative sequence. 
Given the relationship between $f$ and $\phi$, we rewrite \rightInverseDef\ as
$$
  g(z)=z+\sum_{n\geq 2} \phi_{n-1} g(z)g(z/q)\cdots g(z/q^{n-1}) \, .
  \eqno{\equa{gDivergentA}}
$$
As we have seen, \gntn\ and \Segner\ imply that the sequence $(g_n)$ is
nonnegative. Let $R$ be the radius of convergence of $g$. Since  the function
$\phi$ is not degenerate, one of its coefficients, say, $\phi_N$, is positive.
Then \gDivergentA\ yields coefficientwise for formal power series in $x$,
$$
  g(x)\geq \phi_Ng(x)g(x/q)\cdots g(x/q^N)) \, .
  \eqno{\equa{gDivergentB}}
$$
Whenever $x$ is in $[0,R)$, the left hand side of \gDivergentB\
is finite, forcing the right hand side to be finite as well. In other words,
if $x$ is in $[\,0,R)$ then $x/q^N$ is in $[\,0,R)$. Since $q$ is positive
and less than $1$, this forces
$R$ to be either $0$ or $+\infty$. 

Assume that $R$ is infinite. Since $(g_n)$ is nonnegative, $g$ tends to 
infinity as $x$ tends to infinity. Dividing both sides of \gDivergentB\ by
$g(x)$ and taking limit as $x$ tends to infinity yields a contradiction.
Therefore, $R$ is $0$ and $g$ is a divergent series.

When $\bigl((-1)^{n+1}\phi_n\bigl)$ is a nonnegative series, we consider
$g(-x)$ instead of $g(x)$, as in the proof of Theorem \cor.


\finetune{\vfill\eject}

\noindent{\bf References.}
\medskip

{\leftskip=\parindent \parindent=-\parindent
 \par

C.R.\ Adams (1929). On the linear ordinary $q$-difference equations, {\sl Ann.\
Math.}, 30, 195--205.

C.R.\ Adams (1931). Linear $q$-difference equations, {\sl Bull.\ Amer.\ Math.\
Soc.}, 37, 361--400.

G.\ Andrews (1975). Identities in combinatorics. II: A $q$-analog of the
Lagrange inversion theorem, {\sl Proc.\ Amer.\ Math.\ Soc.}, 53, 240--245.

Ph.\ Barbe, W.P.\ McCormick (2012). $q$-Catalan bases and their dual
coefficients, {\tt arXiv:1211.6206}.

Ph.\ Barbe, W.P.\ McCormick (2013). On $q$-algebraic equations and their
power series solutions, {\tt arXiv:1311.5549}.

E.\ Bender (1974). Asymptotic methods in enumeration, {\sl SIAM Rev.}, 16,
485--515.

J.-P.\ B\'ezivin (1992). Sur les \'equations fonctionnelles aux 
$q$-diff\'erences, {\sl Aequationes Math.}, 43, 159--176.

N.H.\ Bingham, C.M.\ Goldie, J.L.\ Teugels (1989). {\sl Regular Variation}, 
2nd ed. Cambridge University Press.

G.\ Birkhoff (1913). The generalized Riemann problem for linear differential
equations and the allied problems for linear difference and $q$-difference
equations, {\sl Proc.\ Amer.\ Acad.}, 49, 521--568.

R.\ Bojanic, E.\ Seneta (1973). Slowly varying functions and asymptotic 
relations, {\sl J.\ Math.\ Anal.\ Appl.}, 303, 302--315.

P.\ Brak, T.\ Prellberg (1995). Critical exponents from nonlinear functional
equations for partially directed cluster models, {\sl J.\ Stat.\ Phys.},
78, 701--730.

J.\ Cano, P.\ Fortuny Ayuso (2012). Power series solutions of nonlinear
$q$-difference equations and the Newton-Puiseux polygon, {\tt arxiv:1209.0295}.

L.\ Carlitz (1972). Sequences, paths, ballot numbers, {\sl Fibonacci Quart.},
10, 531--549.

L.\ Carlitz, J.\ Riordan (1964). Two element lattice permutation numbers and
their $q$-generalization, {\sl Duke J.\ Math.}, 31, 371--388.

R.D.\ Carmichael (1912). The general theory of linear $q$-difference equations,
{\sl Amer.\ J.\ Math.}, 34, 147--168.

B.\ Drake (2009). Limit of areas under lattice paths, {\sl Discrete Math.},
309, 3936--3953.

P.\ Erd\"{o}s, H.\ Pollard, W.\ Feller (1949). A property of power series
with positive coefficients, {\sl Bull.\ Amer.\ Math.\ Soc.}, 55, 201--204.

W.\ Feller (1968). {\sl An Introduction to Probaility Theory and its 
Applications}, Wiley.

Ph.\ Flajolet, R.\ Sedgewick (2009). {\sl Analytic Combinatorics}, Cambridge.

J.\ F\"urlinger, J.\ Hofbauer (1985). $q$-Catalan numbers, {\sl J.\ Comb.\
Th., A}, 248--264.

A.M.\ Garsia (1981). A $q$-analogue of the Lagrange inversion formula,
{\sl Houston J.\ Math.}, 7, 205--237.

A.M.\ Garsia, M.\ Haiman (1996). A remarkable $q,t$-Catalan sequence 
and $q$-Lagrange inversion, {\sl J.\ Algebraic Combin.}, 5, 191--244.

A.\ Garsia, J.\ Lamperti (1962--63). A discrete renewal theorem with infinite
mean, Commentarii Mathematici Helvitici, 37, 221--234.

I.\ Gessel (1980). A noncommutative generalization and $q$-analog of the 
Lagrange inversion formula, {\sl Trans.\ Amer.\ Math.\ Soc.}, 257, 455-482.

G.H.\ Hardy, J.E.\ Littlewood, G.\ P\'olya (1952). {\sl Inequalities}, 
Cambridge University Press.


Ch.\ Krattenthaler (1988). Operator methods and Lagrange inversion: a unified
approach to Lagrange formulas, {\sl Trans.\ Amer.\ Math.\ Soc.}, 305, 431--465.

Chr.\ Mazza, D.\ Piau (2002). Product of correlated symmetric matrices and
$q$-Catalan numbers, {\sl Prob.\ Theor.\ Rel.\ Fields}, 124, 574--594.

I.G.\ Macdonald (1995). {\sl Symmetric Functions and Hall Polynomials}, Oxford
University Press.

J.-P.\ Ramis (1992). About the growth of entire functions solutions of linear
algebraic $q$-difference equations, {\sl Ann.\ Fac.\ Sci.\ Toulouse},  6 (1),
53--94.

J.-P.\ Ramis, J.\ Sauloy, C.\ Zhang (2013). {\sl Local Analytic Classification
of $q$-Difference Equations}, {\tt arXiv:0903.0853}.

J.\ Sauloy (2000). Syst\`emes aux $q$-diff\'erences singuliers r\'eguliers:
classification, matrice de connexion et monodromie, {\sl Ann.\ Inst.\ Fourier}, 50, 1021--1071.

J.\ Sauloy (2003). Galois theory of Fuchsian $q$-difference equations,
{\sl Ann.\ Sci.\ \'Ecole Norm.\ Sup.}, 36, 925--968.

W.J.\ Trjitzinsky (1938). Theory of non-linear $q$-difference systems,
{\sl Ann.\ Math. Pura Appl.}, 17, 59--106.

C.\ Zhang (1998). Sur un th\'eor\`eme de Maillet-Malgrange pour les \'equations
$q$-diff\'erentielles, {\sl Asymptot.\ Anal.}, 17, 309--314.

C.\ Zhang (1999).. D\'eveloppements asymptotiques $q$-Gevrey et s\'eries 
$Gq$-sommables, {\sl Ann.\ Inst.\ Fourier}, 49, 227--261.

C.\ Zhang (2002). A discrete summation for linear $q$-difference equations
with analytic coefficients: general theory and examples, in Braaksma et al.\
ed., {\sl Differential Equations and the Stokes Phenomenon, Proceedings of the 
conference, Groningen, Netherlands, May 28--30, 2001}, World Scientific.

}

\bigskip\bigskip

\setbox1=\vbox{\halign{#\hfil &\hskip 40pt #\hfill\cr
  Ph.\ Barbe                  & W.P.\ McCormick\cr
  90 rue de Vaugirard         & Dept.\ of Statistics \cr
  75006 PARIS                 & University of Georgia \cr
  FRANCE                      & Athens, GA 30602 \cr
  philippe.barbe@math.cnrs.fr & USA \cr
                              & bill@stat.uga.edu \cr}}
\box1

\DoNotPrint{

\bigskip

{\leftskip=\parindent \parindent=-\parindent
 \par
N.H.\ Bingham, R.A.\ Doney (1974). Asymptotic properties of supercritical
branching processes. I: The Galton-Watson process, {\sl Adv.\ Appl.\ Probab.},
6, 711--731.

N.H.\ Bingham, J.\ Hawkes (1983). Some limit theorems for occupation times,
{\sl Probability, Statistics \& Analysis} (eds.\ J.F.C.\ Kigman and G.E.H.\
Reuter), pp.46--62, London Math.\ Soc.\ Lecture Notes, vol.\ 79, CUP.

J.L.\ Teugels (1977). On the rate of convergence of the maximum of a compound
Poisson process, {\sl Bull.\ Soc.\ Math.\ Belg.}, 29, 205--216.

}

\vfill\eject

}

\bye